\documentclass[twocolumn]{autart}
\usepackage{graphics}
\usepackage{epsfig}
\usepackage{times}
\usepackage{amsmath}
\usepackage{amssymb}
\usepackage{graphicx}
\usepackage{epstopdf}
\usepackage{epsfig}
\usepackage{enumerate}
\usepackage[noadjust]{cite}
\usepackage{color}
\usepackage{makecell}
\usepackage{arydshln}
\usepackage{stmaryrd}
\usepackage{stackengine}
\usepackage{threeparttable}
\usepackage{algorithm}
\usepackage{algorithmicx}
\usepackage{xpatch}
\makeatletter
\xpatchcmd\ALG@step{\arabic{ALG@line}}{\fmtlinenumber{ALG@line}}{}{}
\makeatother 
\let\fmtlinenumber\arabic 
\newcommand\mathalph[1]{$\alph{#1}$} 

\newtheorem{Theorem}{Theorem}
\newtheorem{Lemma}{Lemma}
\newtheorem{Definition}{Definition}
\newtheorem{Remark}{Remark}

\newtheorem{Assumption}{Assumption}

\begin{document}
\begin{frontmatter}
\runtitle{Dynamics based Privacy Preservation in Decentralized Optimization}

\title{Dynamics based Privacy Preservation in Decentralized Optimization\thanksref{grants}}

\thanks[grants]{The work was supported in part by the National Science Foundation under Grants ECCS-1912702 and CCF-2106293.}
\author[NPU]{Huan Gao\thanksref{first-author}}\ead{huangao@nwpu.edu.cn}, 
\author[Clemson]{Yongqiang Wang\thanksref{cor-author}}\ead{yongqiw@clemson.edu},
\author[ASU]{Angelia Nedi\'c}\ead{angelia.nedich@asu.edu} 
\thanks[first-author]{Huan Gao was with the Department of Electrical and Computer Engineering, Clemson University, Clemson, SC 29634, USA. He is now with the School of Automation, Northwestern Polytechnical University, Xi'an 710129, China. The work was done when Huan Gao was with Clemson University.}
\thanks[cor-author]{Corresponding author.}
\address[NPU]{School of Automation, Northwestern Polytechnical University, Xi'an 710129, China}
\address[Clemson]{Department of Electrical and Computer Engineering, Clemson University, Clemson, SC 29634, USA}
\address[ASU]{School of Electrical, Computer and Energy Engineering, Arizona State University, Tempe, AZ 85281, USA}


\begin{abstract}
	With decentralized optimization having increased applications in various domains ranging from machine learning, control, sensor networks, to robotics, its privacy is also receiving increased attention. Existing privacy-preserving approaches for decentralized optimization achieve privacy preservation by patching decentralized optimization with information-technology privacy mechanisms such as differential privacy or homomorphic encryption, which either sacrifices optimization accuracy or incurs heavy computation/communication overhead. We propose an inherently privacy-preserving decentralized optimization algorithm by exploiting the robustness of decentralized optimization to uncertainties in optimization dynamics. More specifically, we present a general decentralized optimization framework, based on which we show that privacy can be enabled in decentralized optimization by adding randomness in optimization parameters. We further show that the added randomness has no influence on the accuracy of optimization, and prove that our inherently privacy-preserving algorithm has $R$-linear convergence when the global objective function is smooth and strongly convex. We also rigorously prove that the proposed algorithm can avoid the gradient of a node from being inferable by other nodes. Numerical simulation results confirm the theoretical predictions. 
\end{abstract}
\end{frontmatter}

\section{Introduction}

Decentralized optimization has received increased attention due to its vast applications in online learning \cite{yan2012distributed}, distributed sensing \cite{bazerque2010distributed}, formation control \cite{raffard2004distributed}, source localization \cite{zhang2018distributed}, and power system control \cite{gan2013optimal}. In many of these applications, a network of nodes collectively solve the following problem
\begin{equation}\label{opti_problem}
\begin{aligned}
\min_{\mathbf{x}\in \mathbb{R}^d} \ \ F(\mathbf{x}) \triangleq \sum_{i=1}^{n} f_i(\mathbf{x})
\end{aligned}
\end{equation}
where $f_i : \mathbb{R}^d \rightarrow \mathbb{R}$ is a local objective function accessible only to node $i$.

Over the past decade, a number of gradient-based first order algorithms have been developed to solve the problem. Early results include the decentralized subgradient (DGD) algorithm \cite{nedic2009distributed} which combines average consensus with (sub)gradient descent under diminishing stepsizes. Its convergence rate is $\mathcal{O}({(\ln k)}/{\sqrt{k}})$ for general convex functions and $\mathcal{O}((\ln k)/{k})$ for strongly convex functions, where $k$ is the number of iterations. \cite{yuan2016convergence} shows that the convergence rate of DGD can be improved under a fixed stepsize but at the expense of optimization accuracy. To guarantee both fast convergence and exact solution under fixed stepsizes, many algorithms propose to replace the local gradients in DGD with an auxiliary variable which tracks the global gradient, with typical examples include Aug-DGM \cite{xu2015augmented}, DIGing \cite{qu2017harnessing, nedic2017achieving}, ATC-DIGing \cite{nedic2017geometrically}, AsynDGM \cite{xu2017convergence}, $\mathcal{A}\mathcal{B}$ \cite{xin2018linear}, and Push-Pull \cite{pu2018push, du2018accelerated, zhang2019computational}. These algorithms can achieve $R$-linear convergence.\footnote{Suppose that a sequence $\{\mathbf{x}^k\}$ converges to $\mathbf{x}^*$ in some norm $|| \cdot ||$. Then the convergence is $R$-linear with decay parameter $\rho\in (0, 1)$ if there exists a positive constant $c$ such that $||\mathbf{x}^k - \mathbf{x}^*|| \leq c \rho^k$ holds for any $k$ \cite{nedic2017achieving}.} While most of existing algorithms assume a constant interaction graph, convergence of such gradient-tracking based algorithms on time-varying graphs have also been discussed in \cite{nedic2017achieving}, \cite{xu2017convergence}, and \cite{saadatniaki2020decentralized}. Other relevant algorithms include \cite{shi2015extra, xi2018add, nedic2015distributed, hale2017cloud, fazlyab2018analysis}.

However, none of the aforementioned algorithms consider the privacy of individual nodes, which is unacceptable in many applications. For example, in the rendezvous problem where a group of nodes use decentralized optimization to agree on the optimal assembly position, individual nodes may want to keep their initial positions private, which is important in hostile environments. As indicated in \cite{huang2015differentially}, without protection by an appropriate privacy mechanism, a node's initial position can be easily inferred by an adversary in gradient-decent based rendezvous algorithms. Another example illustrating the importance of privacy in decentralized optimization is collaborative machine learning where gradients/model updates exchanged among participating machines may contain sensitive information such as personal medical record and salary \cite{yan2012distributed}.

To address the pressing need of privacy in decentralized optimization, recently results have emerged on privacy-preserving decentralized optimization. For example, differential-privacy based approaches are proposed in \cite{nozari2016differentially, huang2015differentially, wang2022tailoring}. However, such approaches will unavoidably compromise the accuracy of optimization results. To enable privacy protection with guaranteed optimization accuracy, partially homomorphic encryption based approaches have been proposed in our own prior results \cite{zhang2018enabling, zhang2018admm} as well as others' \cite{lu2018privacy}. However, such approaches will incur heavy computation and communication overhead. \cite{yan2012distributed} and \cite{lou2018privacy} showed that privacy can be obtained by incorporating a projection step or injecting constant uncertainties in stepsizes. However, both approaches have limitations in privacy protection: projection based defense requires individual agents to have a priori knowledge of the optimal solution, whereas constant uncertainties in stepsizes are unable to cover arbitrarily large variations on the gradients. Other approaches include \cite{gade2018private, li2020privacy, wang2022decentralized}. However, they are only applicable to undirected graphs.

Recently, through using random coefficients and/or initial conditions, others as well as our group have proposed several private consensus algorithms \cite{manitara2013privacy, charalambous2019privacy, pilet2019robust, hadjicostis2020privacy, kia2015dynamic, mo2017privacy, gupta2017privacy, he2018distributed, Huan2018CNS, ruan2019secure, ridgley2019simple}. These algorithms can protect the initial value of a node from being inferable by adversaries. Inspired by this line of research, in this paper, we propose to protect the gradient of participating nodes in decentralized optimization by leveraging the robustness of decentralized optimization dynamics. More specifically, by judiciously injecting uncertainties in optimization dynamics, we obfuscate exchanged information without affecting convergence to the exact optimal solution. We rigorously prove that the proposed algorithm can avoid the gradient of a node from being inferable by other nodes. Since protecting the gradient means protecting the values of the gradient function over the entire domain (or protecting both function types and function parameters), it is much more challenging than protecting a single initial value considered in the private consensus problem.

The main contributions of the paper are as follows: 1) We propose a dynamics based privacy protection approach for decentralized optimization that neither sacrifices optimization accuracy nor incurs heavy computation/communication overhead. This is in distinct difference from existing approaches based on differential privacy (which compromise optimization accuracy) and approaches based on homomorphic encryption (which incur heavy computation and communication overhead); 2) Our approach is also different from the multi-party secure computation approach in \cite{he2020secure} which requires each node to communicate with two non-colluding external servers. By introducing randomness into interaction parameters, our approach is implementable in a fully decentralized manner without the assistance of any external servers; 3) We propose a new privacy definition based on the indistinguishability of gradient's arbitrary variations to adversaries from the viewpoint of accessible information, which is stricter than the unobservability/unsolvability based privacy definitions; 4) To facilitate the dynamics based privacy design, we propose a general framework for gradient-tracking based decentralized optimization which includes as special cases many existing decentralized optimization algorithms, such as Aug-DGM \cite{xu2015augmented}, DIGing \cite{qu2017harnessing, nedic2017achieving}, ATC-DIGing \cite{nedic2017geometrically}, AsynDGM \cite{xu2017convergence}, $\mathcal{A}\mathcal{B}$ \cite{xin2018linear}, and Push-Pull \cite{pu2018push, du2018accelerated, zhang2019computational}; 5) We analytically prove that despite the time-varying randomness injected into interaction parameters and the general framework, our proposed approach can still maintain $R$-linear convergence when the global objective function is strongly convex.

The rest of this paper is organized as follows. Sec. 2 formally introduces the problem and reviews some key concepts. A general framework for gradient-tracking based decentralized optimization is introduced in Sec. 3, followed by the inherently privacy-preserving algorithm in Sec. 4. In Sec. 4 we also analyze the convergence of the proposed algorithm and rigorously characterize its privacy-preserving performance. Examples and comparison with existing results are presented in Sec. 5. Finally we conclude this paper in Sec. 6.

{\textbf{Notations}}: $\mathbb{R}$ and $\mathbb{Z}_{\geq 0}$ denote real numbers and nonnegative integers, respectively. $\mathbb{R}^n$ denotes the Euclidean space of dimension $n$, and $\mathbb{R}^{m \times n}$ denotes the set of $m \times n$ matrices with real coefficients. $\mathbf{0}_n \in \mathbb{R}^n$ and $\mathbf{0}_{m \times n} \in \mathbb{R}^{m \times n}$ denote zero vector and $m \times n$ zero matrix, respectively. $\mathbf{1}_n \in \mathbb{R}^n$ denotes the $n \times 1$ all-ones vector, $\mathbf{1}_{m \times n} \in \mathbb{R}^{m \times n}$ denotes the $m \times n$ all-ones matrix, and $\mathbf{I}_n \in \mathbb{R}^{n \times n}$ denotes the $n\times n$ identity matrix. For matrix $\mathbf{A} \in \mathbb{R}^{m \times n}$, its transpose is denoted as $\mathbf{A}^T \in \mathbb{R}^{n \times m}$. The notation $\|\cdot\|$ denotes the Euclidean norm of vectors and the spectral norm of matrices. The notation $\otimes$ represents the Kronecker product.

\section{Problem Formulation}

We characterize the interaction among nodes as a directed graph $\mathcal{G} = (\mathcal{V}, \, \mathcal{E})$, where $\mathcal{V} = \{1, 2, \ldots, n\}$ is the index set of nodes. $\mathcal{E} \subset \mathcal{V} \times \mathcal{V}$ is the set of edges, whose elements are such that an ordered pair $(i, \, j)$ belongs to $\mathcal{E}$ if and only if there exists a directed link from node $j$ to node $i$, i.e., node $j$ can send messages to node $i$. For notational convenience, we assume no self edges, i.e., $(i, \, i)\notin \mathcal{E}$ for all $i\in \mathcal{V}$. The out-neighbor set of node $i$, which represents the set of nodes that can receive messages from node $i$, is denoted as $\mathcal{N}_i^{out}=\{j \in \mathcal{V} \, | \, \forall \, (j, \, i)\in \mathcal{E}\}$. Similarly, the in-neighbor set of node $i$, which represents the set of nodes that can send messages to node $i$, is denoted as $\mathcal{N}_i^{in}=\{j \in \mathcal{V} \, | \, \forall \, (i, \, j)\in \mathcal{E}\}$. From the above definitions, one can obtain that $i \in \mathcal{N}_j^{out}$ and $j \in \mathcal{N}_i^{in}$ are equivalent to each other.

\begin{Assumption}\label{Assumption_directed_graph}
	We assume that the directed graph $\mathcal{G}$ is strongly connected, i.e., for any $i, j \in \mathcal{V}$ with $i\neq j$, there exists at least one directed path from $i$ to $j$ in $\mathcal{G}$, where the directed path respects the direction of the edges.
\end{Assumption}

We use the following standard definitions to characterize objective functions:
\begin{Definition}\label{definition_general_convex}
	A differentiable function $f : \mathbb{R}^d \rightarrow \mathbb{R}$ is convex if $f(\mathbf{x}')\geq f(\mathbf{x}) + \nabla f(\mathbf{x})^T (\mathbf{x}' - \mathbf{x})$ holds for any $\mathbf{x}$ and $\mathbf{x}'$ in $\mathbb{R}^d$.
\end{Definition}

\begin{Definition}\label{definition_strongly_convex}
	A differentiable function $f : \mathbb{R}^d \rightarrow \mathbb{R}$ is $\alpha$-strongly convex with $\alpha > 0$ if $(\nabla f(\mathbf{x})-\nabla f(\mathbf{x}'))^T (\mathbf{x}-\mathbf{x}') \geq \alpha \|{\mathbf{x}-\mathbf{x}'}\|^2$ holds for any $\mathbf{x}$ and $\mathbf{x}'$ in $\mathbb{R}^d$.
\end{Definition}

\begin{Definition}\label{definition_Lipschitz_continuous}
	A differentiable function $f : \mathbb{R}^d \rightarrow \mathbb{R}$ is $\beta$-smooth with $\beta > 0$ if $\|\nabla f(\mathbf{x})-\nabla f(\mathbf{x}')\| \leq \beta \|{\mathbf{x}-\mathbf{x}'}\|$ holds for any $\mathbf{x}$ and $\mathbf{x}'$ in $\mathbb{R}^d$.
\end{Definition}

\begin{Assumption}\label{Assumption_objective_functions}
	We assume that each local objective function $f_i$ is differentiable, convex, and $\beta_i$-smooth. We also assume that the global objective function $F=\sum_{i=1}^{n}f_i(\mathbf{x})$ is $\alpha_F$-strongly convex.
\end{Assumption}

Under Definitions \ref{definition_general_convex}-\ref{definition_Lipschitz_continuous} and Assumption \ref{Assumption_objective_functions}, one can verify that the global objective function $F$ is $\beta_F$-smooth with $\beta_F=\sum_{i=1}^{n}\beta_i$, and $\alpha_F$ is always less than or equal to $\beta_F$. Moreover, problem (\ref{opti_problem}) has a unique optimal solution, which is denoted by $\mathbf{x}^*$.

\section{A General Decentralized Optimization Framework}

\subsection{A New Framework for Decentralized Optimization}

As discussed in Sec. 1, most existing decentralized optimization algorithms do not consider the privacy of individual nodes. In this paper, we propose to enable privacy preservation in decentralized optimization by exploiting the robustness of decentralized optimization dynamics. To this end, we first propose a new framework for gradient-tracking based decentralized optimization in Algorithm 1.

\begin{algorithm}
	\caption{A new decentralized optimization framework}
	Each node $i$ initializes $\mathbf{x}_i^{0}$ randomly in $\mathbb{R}^d$ and sets $\mathbf{y}_i^{0} = \nabla f_i(\mathbf{x}_i^{0})$. At iteration $k$:
	\begin{algorithmic}[1]
		\let\fmtlinenumber\mathalph
		\State Node $i$ computes and sends $\mathbf{x}_i^{k}$ as well as $\mathbf{\Lambda}_i^{k} \mathbf{y}_i^{k}$ to its out-neighbors $l \in \mathcal{N}_i^{out}$, where $\mathbf{\Lambda}_i^{k}$ is a diagonal matrix denoting the stepsize.
		\State After receiving $\mathbf{x}_j^{k}$ and $\mathbf{\Lambda}_j^{k} \mathbf{y}_j^{k}$ from its in-neighbors $j \in \mathcal{N}_i^{in}$, node $i$ updates $\mathbf{x}_i$ as:
		\begin{equation}\label{update-rule-01}
			\begin{aligned}
				\mathbf{x}_i^{k+1} = {\textstyle \sum_{j\in \mathcal{N}_i^{in} \cup \{i\}}} \big({\mathbf{R}_{ij}^{k} \mathbf{x}_j^{k} - \mathbf{A}_{ij}^{k} \mathbf{\Lambda}_j^{k} \mathbf{y}_j^{k}}\big)
			\end{aligned}
		\end{equation}
		where $\mathbf{R}_{ij}^{k}$ and $\mathbf{A}_{ij}^{k}$ are coupling weight matrices.
		\State After updating $\mathbf{x}_i$, node $i$ computes and sends $\mathbf{C}_{li}^{k} \mathbf{y}_i^{k} + \mathbf{B}_{li}^{k} \big( \nabla f_i(\mathbf{x}_i^{k+1}) - \nabla f_i(\mathbf{x}_i^{k})\big)$ to its out-neighbors $l \in \mathcal{N}_i^{out}$, where $\mathbf{C}_{li}^{k}$ and $\mathbf{B}_{li}^{k}$ are coupling weight matrices.
		\State After receiving $\mathbf{C}_{ij}^{k} \mathbf{y}_j^{k} + \mathbf{B}_{ij}^{k} \big( \nabla f_j(\mathbf{x}_j^{k+1}) - \nabla f_j(\mathbf{x}_j^{k})\big)$ from its in-neighbors $j \in \mathcal{N}_i^{in}$, node $i$ updates $\mathbf{y}_i$ as:
		\begin{equation}\label{update-rule-02}
			\begin{aligned}
				\mathbf{y}_i^{k+1} = & {\textstyle \sum_{j\in \mathcal{N}_i^{in} \cup \{i\}}} \big( \mathbf{C}_{ij}^{k} \mathbf{y}_j^{k} \\
				& \quad + \mathbf{B}_{ij}^{k} \big( \nabla f_j(\mathbf{x}_j^{k+1}) - \nabla f_j(\mathbf{x}_j^{k})\big) \big)
			\end{aligned}
		\end{equation}
	\end{algorithmic}
\end{algorithm}

Note that by setting $\mathbf{R}_{ij}^{k} = \mathbf{A}_{ij}^{k} = \mathbf{0}_{d\times d}$ for $j\notin \mathcal{N}_i^{in} \cup \{i\}$ and $\mathbf{C}_{ji}^{k} = \mathbf{B}_{ji}^{k} = \mathbf{0}_{d\times d}$ for $j\notin \mathcal{N}_i^{out} \cup \{i\}$, the update rules in (\ref{update-rule-01})-(\ref{update-rule-02}) can be rewritten as
\begin{equation}\label{update-rule-2}
	\begin{aligned}
		\mathbf{x}_i^{k+1} & = {\textstyle \sum_{j=1}^n} \big({\mathbf{R}_{ij}^{k} \mathbf{x}_j^{k} - \mathbf{A}_{ij}^{k} \mathbf{\Lambda}_j^{k} \mathbf{y}_j^{k}}\big)\\
		\mathbf{y}_i^{k+1} & = {\textstyle \sum_{j=1}^n} \big( \mathbf{C}_{ij}^{k} \mathbf{y}_j^{k} + \mathbf{B}_{ij}^{k} \big( \nabla f_j^{k+1} - \nabla f_j^{k}\big) \big)
	\end{aligned}
\end{equation}
where $\nabla f_j^{k} = \nabla f_j(\mathbf{x}_j^{k})$ is the gradient of $f_j$ evaluated at $\mathbf{x}_j^{k}$, which is typically nonlinear in nature. (\ref{update-rule-2}) can be rewritten into a more compact matrix form
\begin{equation}\label{system_matrix}
\begin{aligned}
\mathbf{x}^{k+1} & = \mathbf{R}^{k} \mathbf{x}^{k} - \mathbf{A}^{k} \mathbf{\Lambda}^{k} \mathbf{y}^{k}\\
\mathbf{y}^{k+1} & = \mathbf{C}^{k} \mathbf{y}^{k} + \mathbf{B}^{k} \big( \nabla f^{k+1} - \nabla f^{k}\big)
\end{aligned}
\end{equation}
where $\mathbf{x}^{k} =[(\mathbf{x}_1^{k})^T \, \cdots \, (\mathbf{x}_n^{k})^T]^T$, $\mathbf{y}^{k} =[(\mathbf{y}_1^{k})^T \, \cdots \, (\mathbf{y}_n^{k})^T]^T$, $\nabla f^{k} =[(\nabla f_1^{k})^T \, \cdots \, (\nabla f_n^{k})^T]^T$, and $\mathbf{R}^{k}$, $\mathbf{A}^{k}$, $\mathbf{C}^{k}$, and $\mathbf{B}^{k}$ are block matrices with the $(ij)$-th block entry being $\mathbf{R}_{ij}^{k}$, $\mathbf{A}_{ij}^{k}$, $\mathbf{C}_{ij}^{k}$, and $\mathbf{B}_{ij}^{k}$, respectively. $\mathbf{\Lambda}^{k}$ is a block diagonal matrix with the $i$-th diagonal block being $\mathbf{\Lambda}_i^{k}$.

\subsection{Relationship with Existing Decentralized Optimization Algorithms}

The proposed decentralized optimization framework is very general and includes as special cases many popular decentralized optimization algorithms, such as Aug-DGM \cite{xu2015augmented}, DIGing \cite{qu2017harnessing, nedic2017achieving}, ATC-DIGing \cite{nedic2017geometrically}, AsynDGM \cite{xu2017convergence}, $\mathcal{A}\mathcal{B}$ \cite{xin2018linear}, and Push-Pull \cite{pu2018push, du2018accelerated, zhang2019computational}. Table I summarizes the particular selections of parameters $\mathbf{R}^{k}$, $\mathbf{A}^{k}$, $\mathbf{B}^{k}$, $\mathbf{C}^{k}$, and $\mathbf{\Lambda}^{k}$ in (\ref{system_matrix}) to obtain some commonly used decentralized optimization algorithms. 

\begin{table}[h!]
	\centering{Table I. Particular selections of parameters in the \\ proposed framework to obtain some existing algorithms.}
	\begin{center}
		\begin{threeparttable}	
			\renewcommand\arraystretch{1.6}
			\begin{tabular}{|c|c|c|c|c|c|}
				\hline & $\mathbf{R}^{k}$ & $\mathbf{A}^{k}$ & $\mathbf{C}^{k}$ & $\mathbf{B}^{k}$ & $\mathbf{\Lambda}^{k}$ \\ \hline
				Aug-DGM \cite{xu2015augmented} & $\mathbf{W}$ & $\mathbf{W}$ & $\mathbf{W}$ & $\mathbf{W}$ & $\mathbf{\Lambda}$ \\ \hline
				DIGing \cite{qu2017harnessing} & $\mathbf{W}$ & $\mathbf{I}$ & $\mathbf{W}$ & $\mathbf{I}$ & $\lambda\mathbf{I}$ \\ \hline
				DIGing \cite{nedic2017achieving} & $\mathbf{W}^{k}$ & $\mathbf{I}$ & $\mathbf{W}^{k}$ & $\mathbf{I}$ & $\lambda\mathbf{I}$ \\ \hline
				ATC-DIGing \cite{nedic2017geometrically} & $\mathbf{W}$ & $\mathbf{W}$ & $\mathbf{W}$ & $\mathbf{W}$ & $\mathbf{\Lambda}$ \\ \hline
				AsynDGM \cite{xu2017convergence} & $\mathbf{W}^{k}$ & $\mathbf{W}^{k}$ & $\mathbf{W}^{k}$ & $\mathbf{I}$ & $\mathbf{\Lambda}$ \\ \hline
				$\mathcal{A}\mathcal{B}$ \cite{xin2018linear} & $\mathbf{R}$ & $\mathbf{I}$ & $\mathbf{C}$ & $\mathbf{C}$ & $\lambda\mathbf{I}$ \\ \hline
				Push-Pull \cite{pu2018push} & $\mathbf{R}$ & $\mathbf{R}$ & $\mathbf{C}$ & $\mathbf{C}$ & $\lambda\mathbf{I}$ \\ \hline
				Push-Pull \cite{du2018accelerated} & $\mathbf{R}$ & $\mathbf{I}$ & $\mathbf{C}$ & $\mathbf{I}$ & $\lambda\mathbf{I}$ \\ \hline
				Push-Pull \cite{zhang2019computational} & $\mathbf{R}$ & $\mathbf{R}$ & $\mathbf{C}$ & $\mathbf{I}$ & $\lambda\mathbf{I}$ \\ \hline	
			\end{tabular}
		\end{threeparttable}
	\end{center}
	\begin{tablenotes} 	
		\item [] $\mathbf{W}$ and $\mathbf{W}^{k}$ are doubly stochastic, $\mathbf{R}$ is time-invariant and row-stochastic, $\mathbf{C}$ is time-invariant and column-stochastic, and $\lambda\mathbf{I}$ and $\mathbf{\Lambda}$ represent homogeneous and heterogeneous stepsize matrices, respectively. Note that Kronecker product with $\mathbf{I}_d$ is needed for these matrices when $d>1$.
	\end{tablenotes} 
\end{table}

In fact, by setting $\mathbf{R}^{k}$, $\mathbf{A}^{k}$, $\mathbf{C}^{k}$, $\mathbf{B}^{k}$, and $\mathbf{\Lambda}^{k}$ to $\mathbf{R}$, $\mathbf{I}$, $\mathbf{C}$, $\mathbf{I}$, and $\mathbf{\Lambda}$, respectively, our proposed algorithm in (\ref{system_matrix}) can be rewritten as $\mathbf{x}^{k+1} = (\mathbf{R}+\mathbf{C}) \mathbf{x}^{k} -\mathbf{C} \mathbf{R}\mathbf{x}^{k-1} - \mathbf{\Lambda} \big( \nabla f^{k} - \nabla f^{k-1}\big)$ which becomes EXTRA \cite{shi2015extra}.

As indicated in Table I, our proposed framework reduces to existing algorithms when the parameters are selected appropriately. Moreover, it can also give rise to new algorithms with special properties. In particular, in what follows we show that it results in new algorithms having inherent privacy-preserving capabilities.

\section{Privacy-preserving Decentralized Optimization}

Based on the general decentralized optimization framework proposed in Sec. 3, we can enable privacy-preservation in decentralized optimization without compromising optimization accuracy.

\subsection{Privacy-preserving Design}

To enable privacy, we propose to add randomness in stepsize $\mathbf{\Lambda}^{k}$ and coupling weights $\mathbf{R}^{k}$, $\mathbf{A}^{k}$, $\mathbf{C}^{k}$, and $\mathbf{B}^{k}$ for iterations $k< K$, where $K$ is a positive integer. The detailed parameter design for each node $i$ is given in Table II. More specifically, for iterations $k< K$, each node $i$ selects the coupling weights $\mathbf{C}_{ji}^k = \text{diag} \left\{ c_{ji}^k(1), \ldots, c_{ji}^k(d) \right\}$ and $\mathbf{B}_{ji}^k = \text{diag} \left\{ b_{ji}^k(1), \ldots, b_{ji}^k(d) \right\}$ for $j \in \mathcal{N}_i^{out}$ following any chosen random distributions with support $\mathbb{R}$ such as Gaussian or Laplace distribution (so the coupling weights can be negative, positive, or zero). To guarantee the column-stochastic property of both $\mathbf{C}^k$ and $\mathbf{B}^k$, each node $i$ sets $\mathbf{C}_{ii}^k$ and $\mathbf{B}_{ii}^k$ as $\mathbf{C}_{ii}^k=\mathbf{I}_d-\sum_{j\in \mathcal{N}_i^{out}} \mathbf{C}_{ji}^k$ and $\mathbf{B}_{ii}^k=\mathbf{I}_d-\sum_{j\in \mathcal{N}_i^{out}} \mathbf{B}_{ji}^k$, respectively. In this way, the column-stochastic property of both $\mathbf{C}^k$ and $\mathbf{B}^k$ for iterations $k <K$ is guaranteed in a fully decentralized manner. Since we do not require $\mathbf{R}^{k}$ or $\mathbf{A}^{k}$ to be row-stochastic for iterations $k<K$, each node $i$ can select $\mathbf{R}_{ij}^k$ and $\mathbf{A}_{ij}^k$ for $j \in \mathcal{N}_i^{in}\cup\{i\}$ following any random distributions.

\begin{table*}
	\centering{Table II. Parameter design for each node $i$ in our proposed framework.}
	\begin{center}
		\begin{threeparttable}	
			\renewcommand\arraystretch{2.5}
			\begin{tabular}{|c|c |c|}
				\hline & Iterations $k<K \, ^\dag$ & Iterations $k\geq K$ \\
				\hline \makecell{$\mathbf{\Lambda}_i^{k}$} 
				& \makecell{ $\mathbf{\Lambda}_i^{k}= \text{diag}\{\lambda_i^{k}(1),\cdots, \lambda_i^{k}(d)\}$, where $\lambda_i^{k}(1),\ldots, \lambda_i^{k}(d)$ \\ are chosen following selected distributions with support $\mathbb{R}$}
				& \makecell{$\mathbf{\Lambda}_i^{k} = \lambda \mathbf{I}_d$, $\lambda >0$} \\
				\hline \makecell{$\mathbf{R}_{ij}^{k}$} 
				& \makecell{$\mathbf{R}_{ij}^{k} = \text{diag}\{r_{ij}^{k}(1),\cdots, r_{ij}^{k}(d)\}$, where $r_{ij}^{k}(l)$ are chosen from $\mathbb{R}$ for \\ $j\in \mathcal{N}_i^{in} \cup \{i\}$ and $l=1,\ldots,d$ following selected distributions with support $\mathbb{R}$} 
				& \makecell{$\mathbf{R}_{ij}^{k} = r_{ij}^{k} \mathbf{I}_d$, where $r_{ij}^{k}$ are selected from $[\eta, \, 1]$ \\ for $j\in \mathcal{N}_i^{in} \cup \{i\}$ subject to $\sum_{j=1}^n r_{ij}^{k} =1$} \\	
				\hline \makecell{$\mathbf{A}_{ij}^{k}$} 
				& \makecell{$\mathbf{A}_{ij}^{k} = \text{diag}\{a_{ij}^{k}(1),\cdots, a_{ij}^{k}(d)\}$, where $a_{ij}^{k}(l)$ are chosen from $\mathbb{R}$ for \\ $j\in \mathcal{N}_i^{in} \cup \{i\}$ and $l=1,\ldots,d$ following selected distributions with support $\mathbb{R}$}
				& \makecell{$\mathbf{A}_{ij}^{k} = a_{ij}^{k} \mathbf{I}_d$, where $a_{ij}^{k}$ are selected from $[\eta, \, 1]$ \\ for $j\in \mathcal{N}_i^{in} \cup \{i\}$ subject to $\sum_{j=1}^n a_{ij}^{k} =1$} \\
				\hline \makecell{$\mathbf{C}_{ji}^{k}$}
				& \makecell{$\mathbf{C}_{ji}^{k}= \text{diag}\{c_{ji}^{k}(1),\cdots, c_{ji}^{k}(d)\}$, where $c_{ji}^{k}(l)$ are chosen from $\mathbb{R}$ for \\ $j \in \mathcal{N}_i^{out}$ and $l=1,\ldots,d$ following selected distributions with support $\mathbb{R}$, \\ and $\mathbf{C}_{ii}^k$ is set as $\mathbf{C}_{ii}^k=\mathbf{I}_d-\sum_{j\in \mathcal{N}_i^{out}} \mathbf{C}_{ji}^k$}
				& \makecell{$\mathbf{C}_{ji}^{k}= c_{ji}^{k} \mathbf{I}_d$, where $c_{ji}^{k}$ are selected from $[\eta, \, 1]$ \\ for $j \in \mathcal{N}_i^{out} \cup \{i\}$ subject to $\sum_{j=1}^n c_{ji}^{k} =1$} \\
				\hline \makecell{$\mathbf{B}_{ji}^{k}$}
				& \makecell{$\mathbf{B}_{ji}^{k}= \text{diag}\{b_{ji}^{k}(1),\cdots, b_{ji}^{k}(d)\}$, where $b_{ji}^{k}(l)$ are chosen from $\mathbb{R}$ for \\ $j \in \mathcal{N}_i^{out}$ and $l=1,\ldots,d$ following selected distributions with support $\mathbb{R}$, \\ and $\mathbf{B}_{ii}^k$ is set as $\mathbf{B}_{ii}^k=\mathbf{I}_d-\sum_{j\in \mathcal{N}_i^{out}} \mathbf{B}_{ji}^k$}
				& \makecell{$\mathbf{B}_{ji}^{k} = b_{ji}^{k} \mathbf{I}_d$ \\ select $b_{ji}^{k}$ from $[\eta, \, 1]$ for $j \in \mathcal{N}_i^{out} \cup \{i\}$ \\ subject to $\sum_{j=1}^n b_{ji}^{k}=1$} \\
				\hline 		
			\end{tabular}
		\end{threeparttable}
	\end{center}
	\begin{tablenotes}
		\item[] $^\dag$ $K$ can be any positive integer. Its influence will be discussed in detail in Remark \ref{K_convergence} and Remark \ref{K_privacy}.
	\end{tablenotes}
\end{table*}

Furthermore, as indicated in Table II, the column stochastic property of $\mathbf{C}^k$ and $\mathbf{B}^k$ and the row stochastic property of $\mathbf{R}^k$ and $\mathbf{A}^k$ are required for iterations $k\geq K$. Let us take $\mathbf{C}^k$ as an example to show how to ensure these properties in a fully decentralized manner. To meet the requirements of Table II, each node $i$ selects a set of real values $\left\{ c_{ji}^k\in [\eta, \, 1] \, \big| \, j\in \mathcal{N}_i^{out} \cup \{i\} \right\}$ with their sum equal to one, then sets $\mathbf{C}_{ji}^k$ as $\mathbf{C}_{ji}^k = c_{ji}^k \mathbf{I}_d$ for $j\in \mathcal{N}_i^{out} \cup \{i\}$. Note that there are many ways to obtain such a set of real values with sum equal to one. For example, node $i$ first randomly selects a set of real values $\left\{ p_{ji}^k \, \big| \, j\in \mathcal{N}_i^{out} \cup \{i\} \right\}$ from $[0, \, 1]$, then it determines $c_{ji}^k$ by normalizing these values $\{p_{ji}^k\}$ via $c_{ji}^k = \frac{[1- (|\mathcal{N}_i^{out}|+1) \eta] [(1-\eta)p_{ji}^k + \eta]}{(1-\eta)\sum_{l\in \mathcal{N}_i^{out} \cup \{i\}}p_{li}^k + (|\mathcal{N}_i^{out}|+1)\eta} + \eta$ for $j\in \mathcal{N}_i^{out} \cup \{i\}$. One can verify $\sum_{j\in \mathcal{N}_i^{out}\cup \{i\}}c_{ji}^k =1$ and $c_{ji}^k \in [\eta, \, 1]$ for $j\in \mathcal{N}_i^{out} \cup \{i\}$, and hence the column stochastic property of $\mathbf{C}^k$ is guaranteed in a fully decentralized manner.

\begin{Remark}
	We allow each agent $i$ to randomly choose its associated coupling weights and stepsizes for iterations $k<K$ following any distributions. These distributions can be any continuous probability distribution with support $\mathbb{R}$, e.g., Gaussian distribution, Laplace distribution, etc.
\end{Remark}


\subsection{Convergence Analysis}

In this subsection, inspired by the work in \cite{saadatniaki2020decentralized}, we prove that our privacy-preserving algorithm has $R$-linear convergence. Since the parameters in iterations $0 \leq k \leq K-1$ vary randomly with time, we first analyze their influence on states $\mathbf{x}_i^K$ and $\mathbf{y}_i^K$, based on which we can characterize the evolution of state variables after iteration $K-1$. From the parameter design in Table II, we have $\sum_{i=1}^{n}\mathbf{C}_{ij}^{k} =\mathbf{I}_{d}$ and $\sum_{i=1}^{n}\mathbf{B}_{ij}^{k} =\mathbf{I}_{d}$, which, in combination with (\ref{update-rule-2}) leads to
\begin{equation}\label{update-rule-3}
\begin{aligned}
\big({\textstyle \sum_{i=1}^{n}} \mathbf{y}_i^{k+1} \big)^T = \big( {\textstyle \sum_{i=1}^n} (\mathbf{y}_i^{k} + \nabla f_i^{k+1} - \nabla f_i^{k}) \big)^T
\end{aligned}
\end{equation}
for $k \leq K-1$, where we used the property $\big(\mathbf{C}_{ij}^{k}\big)^T = \mathbf{C}_{ij}^{k}$ and $\big(\mathbf{B}_{ij}^{k}\big)^T = \mathbf{B}_{ij}^{k}$ due to their diagonal structure. We can further rewrite (\ref{update-rule-3}) as $\sum_{i=1}^{n} \mathbf{y}_i^{k+1} -\sum_{i=1}^n \nabla f_i^{k+1} = \sum_{i=1}^n \mathbf{y}_i^{k} - \sum_{i=1}^n \nabla f_i^{k}$. Given $\mathbf{y}_i^{0}=\nabla f_i^{0}$, one can obtain
\begin{equation}\label{update-rule-5}
\begin{aligned}
{\textstyle \sum_{i=1}^{n}} \big(\mathbf{y}_i^{K} - \nabla f_i^{K}\big) = \cdots = {\textstyle \sum_{i=1}^n} \big(\mathbf{y}_i^{0} - \nabla f_i^{0}\big) = \mathbf{0}_d
\end{aligned}
\end{equation}

The equality constraint in (\ref{update-rule-5}) reflects the influence of random time-varying parameters in iterations $0 \leq k \leq K-1$ on the dynamics. Next, under this constraint, we study the system dynamics after iteration $K-1$.

For our parameter design in Table II, we have $\mathbf{\Lambda}_i^{k} = \lambda \mathbf{I}_d$, $\mathbf{R}_{ij}^{k}=r_{ij}^{k} \mathbf{I}_d$, $\mathbf{A}_{ij}^{k}=a_{ij}^{k} \mathbf{I}_d$, $\mathbf{C}_{ij}^{k}=c_{ij}^{k} \mathbf{I}_d$, and $\mathbf{B}_{ij}^{k}=b_{ij}^{k} \mathbf{I}_d$ for $k \geq K$. Constructing $n\times n$ matrices $\mathbf{\bar{R}}^{k}$, $\mathbf{\bar{A}}^{k}$, $\mathbf{\bar{C}}^{k}$, and $\mathbf{\bar{B}}^{k}$ with the $(ij)$-th elements being $r_{ij}^{k}$, $a_{ij}^{k}$, $c_{ij}^{k}$, and $b_{ij}^{k}$, respectively, we have $\mathbf{R}^{k}=\mathbf{\bar{R}}^{k} \otimes \mathbf{I}_d$, $\mathbf{A}^{k}=\mathbf{\bar{A}}^{k} \otimes \mathbf{I}_d$, $\mathbf{C}^{k}=\mathbf{\bar{C}}^{k} \otimes \mathbf{I}_d$, and $\mathbf{B}^{k}=\mathbf{\bar{B}}^{k} \otimes \mathbf{I}_d$ for $k \geq K$. Then for $k \geq K$, we can rewrite the system dynamics (\ref{system_matrix}) as
\begin{equation}\label{system-matrix}
\begin{aligned}
\mathbf{x}^{k+1} & = (\mathbf{\bar{R}}^{k} \otimes \mathbf{I}_d) \mathbf{x}^{k} - \lambda (\mathbf{\bar{A}}^{k} \otimes \mathbf{I}_d) \mathbf{y}^{k} \\
\mathbf{y}^{k+1} & =(\mathbf{\bar{C}}^{k} \otimes \mathbf{I}_d) \mathbf{y}^{k} + (\mathbf{\bar{B}}^{k} \otimes \mathbf{I}_d) \big( \nabla f^{k+1} - \nabla f^{k}\big)
\end{aligned}
\end{equation}

To facilitate convergence analysis, we introduce a state transformation, $\mathbf{s}^{k}= ((\mathbf{\bar{V}}^{k})^{-1}\otimes \mathbf{I}_d) \mathbf{y}^{k}$, where $\mathbf{\bar{V}}^{k}=\text{diag}(\mathbf{v}^k)$ with the evolution of $\mathbf{v}^k$ governed by
\begin{equation}\label{Eqn_vector_v}
\mathbf{v}^{k+1}=\mathbf{\bar{C}}^{k} \mathbf{v}^{k}
\end{equation}
for $k \geq K$. The value of $\mathbf{v}^{k}$ at iteration $K$ is set as $\frac{1}{n}\mathbf{1}_n$, i.e., $\mathbf{v}^{K} = \frac{1}{n} \mathbf{1}_n$. Then (\ref{system-matrix}) can be rewritten as
\begin{equation}\label{Eqn_system_matrix_form_3}
	\begin{aligned}
		\mathbf{x}^{k+1} & = \mathbf{R}^{k} \mathbf{x}^{k} - \lambda (\mathbf{\bar{A}}^{k}\mathbf{\bar{V}}^{k} \otimes \mathbf{I}_d) \mathbf{s}^{k} \\
		\mathbf{s}^{k+1} & = \mathbf{P}^{k} \mathbf{s}^{k} + ((\mathbf{\bar{V}}^{k+1})^{-1} \mathbf{\bar{B}}^{k} \otimes \mathbf{I}_d) \big( \nabla f^{k+1} - \nabla f^{k}\big)
	\end{aligned}
\end{equation}
for $k \geq K$ where $\mathbf{P}^{k} = \mathbf{\bar{P}}^{k} \otimes \mathbf{I}_d$ and $\mathbf{\bar{P}}^{k} = (\mathbf{\bar{V}}^{k+1})^{-1} \mathbf{\bar{C}}^{k} \mathbf{\bar{V}}^{k}$. As proven in Lemma \ref{Lemma_P} in Appendix A, $\mathbf{\bar{P}}^{k}$ is row-stochastic and the sequence $\{\mathbf{\bar{P}}^{k}\}$ is ergodic for $k \geq K$. Furthermore, the sequence of $\{\mathbf{v}^k\}$ for $k\geq K$ is an absolute probability sequence for $\{\mathbf{\bar{P}}^{k}\}$ (see Appendix A for the definitions of ergodic sequence and absolute probability sequence.) We define $\mathbf{\bar{x}}_{\text{w}}^k$, $\mathbf{\tilde{x}}_{\text{w}}^k$, $\mathbf{r}^k$, and $\mathbf{\tilde{s}}_{\text{w}}^k$ as
\begin{equation*}
	\begin{aligned}
		\mathbf{\bar{x}}_{\text{w}}^k & =((\boldsymbol{\phi}^k)^T\otimes \mathbf{I}_d) \mathbf{x}^k, \quad \
		\mathbf{r}^k =\mathbf{1}_n\otimes \mathbf{\bar{x}}_{\text{w}}^k - \mathbf{1}_n\otimes \mathbf{x}^* \\
		\mathbf{\tilde{x}}_{\text{w}}^k & = \mathbf{x}^k - \mathbf{1}_n\otimes \mathbf{\bar{x}}_{\text{w}}^k, \qquad
		\mathbf{\tilde{s}}_{\text{w}}^k = \mathbf{s}^k - (\mathbf{1}_n(\mathbf{v}^k)^T \otimes \mathbf{I}_d) \mathbf{s}^k
	\end{aligned}
\end{equation*}
where for $k\geq K$, $\{\boldsymbol{\phi}^k\}$ is an absolute probability sequence for the ergodic sequence of row-stochastic matrices $\{\mathbf{\bar{R}}^k\}$. Before presenting our main result, we introduce some auxiliary lemmas.

\begin{Lemma}\label{Lemma_x}
	Under Assumptions \ref{Assumption_directed_graph} and \ref{Assumption_objective_functions}, and the parameter design in Table II, the following inequality holds for $k\geq K+ \bar{N}-1$:
	\begin{equation}\label{Eqn_lemma_x}
		\begin{aligned}
			 & \|\mathbf{\tilde{x}}_{\text{w}}^{k+1}\| \\
		& \leq \big(r_R + \lambda Q_R n \sqrt{n} \bar{\beta}\big) \|\mathbf{\tilde{x}}_{\text{w}}^{k-\bar{N}+1}\| + \lambda Q_R n \sqrt{n} \bar{\beta} \sum_{l=0}^{\bar{N}-2}	\|\mathbf{\tilde{x}}_{\text{w}}^{k-l}\| \\
		& \quad + \lambda Q_R n \sqrt{n} \bar{\beta} \|\mathbf{r}^{k-\bar{N}+1}\| + \lambda Q_R n \sqrt{n} \bar{\beta} \sum_{l=0}^{\bar{N}-2} \|\mathbf{r}^{k-l}\| \\
		& \quad + \lambda Q_R \sqrt{n} \|\mathbf{\tilde{s}}_{\text{w}}^{k-\bar{N}+1} \| + \lambda Q_R \sqrt{n} \sum_{l=0}^{\bar{N}-2} \|\mathbf{\tilde{s}}_{\text{w}}^{k-l} \|
		\end{aligned}
	\end{equation}
	where $\bar{\beta}=\max\{\beta_1,\cdots, \beta_n\}$, $Q_R=2n(1+\eta^{-(n-1)})/(1-\eta^{n-1})$, $r_{R}=Q_R (1-\eta^{n-1})^{\frac{N_R-1}{n-1}}$, and $\bar{N}=\max\{N_R, N_P\}$ (see Lemmas \ref{Lemma_P} and \ref{Lemma_R} in Appendix A for more details on parameters $N_R$, $N_P$, $r_R$, and $Q_R$.) 
\end{Lemma}

{\it Proof}: The proof is given in Appendix B. \hfill{$\blacksquare$}

\begin{Lemma}\label{Lemma_r}
	Under Assumptions \ref{Assumption_directed_graph} and \ref{Assumption_objective_functions}, and the parameter design in Table II, the following inequality holds for $k\geq K$:
	\begin{equation}\label{Eqn_lemma_r}
		\begin{aligned}
			\|\mathbf{r}^{k+1}\| \leq \lambda n \bar{\beta} \|\mathbf{\tilde{x}}_{\text{w}}^{k}\| + \big(1-\lambda n^{-1}\eta^{n-1} \alpha_F \big) \| \mathbf{r}^k\| + \lambda n \| \mathbf{\tilde{s}}_{\text{w}}^{k}\|
		\end{aligned}
	\end{equation}
	if the stepsize $\lambda$ satisfies $\lambda \leq {1}/\beta_F$.
\end{Lemma}

{\it Proof}: The proof is given in Appendix C. \hfill{$\blacksquare$}

\begin{Lemma}\label{Lemma_s}
	Under Assumptions \ref{Assumption_directed_graph} and \ref{Assumption_objective_functions}, and the parameter design in Table II, the following inequality holds for $k\geq K+ \bar{N}-1$:
	\begin{equation}\label{Eqn_lemma_s}
		\begin{aligned}
			& \|\mathbf{\tilde{s}}_{\text{w}}^{k+1}\| \\
			& \leq \big(\frac{2n^2\bar{\beta}Q_P}{\eta^{n-1}} + \lambda \frac{n^3\bar{\beta}^2Q_P}{\eta^{n-1}}\big) \big(\|\mathbf{\tilde{x}}_{\text{w}}^{k-\bar{N}+1}\| + \sum_{l=0}^{\bar{N}-2}	\|\mathbf{\tilde{x}}_{\text{w}}^{k-l}\|\big) \\
			& \quad + \lambda \frac{n^3\bar{\beta}^2Q_P}{\eta^{n-1}} \big(\|\mathbf{r}^{k-\bar{N}+1}\| + \sum_{l=0}^{\bar{N}-2} \|\mathbf{r}^{k-l}\|\big) \\
			& \quad + \big(r_P + \lambda \frac{n^2\bar{\beta}Q_P}{\eta^{n-1}}\big) \|\mathbf{\tilde{s}}_{\text{w}}^{k-\bar{N}+1} \| + \lambda \frac{n^2\bar{\beta}Q_P}{\eta^{n-1}} \sum_{l=0}^{\bar{N}-2} \|\mathbf{\tilde{s}}_{\text{w}}^{k-l} \|
		\end{aligned}
	\end{equation}
	where $Q_P=2n\big(1+(n\eta^{-n})^{n-1}\big)/\big(1-(n^{-1}\eta^n)^{n-1}\big)$ and $r_{P}= Q_P \big(1-(n^{-1}\eta^n)^{n-1}\big)^{\frac{N_P-1}{n-1}}$ (see Lemma \ref{Lemma_P} in Appendix A for more details on parameters $Q_P$ and $r_P$.)
\end{Lemma}

{\it Proof}: The proof is given in Appendix D. \hfill{$\blacksquare$}

Now we are in position to present our main convergence result.
\begin{Theorem}\label{theorem_convergence}
	Under Assumptions \ref{Assumption_directed_graph} and \ref{Assumption_objective_functions}, and the parameter design in Table II, our algorithm has $R$-linear convergence when the stepsize parameter $\lambda$ is sufficiently small.
\end{Theorem}

{\it Proof}: Denoting $\boldsymbol{\xi}^k$ as $\boldsymbol{\xi}^k = \big[\|\mathbf{\tilde{x}}_{\text{w}}^{k}\|, \, \|\mathbf{r}^{k}\|, \, \|\mathbf{\tilde{s}}_{\text{w}}^{k}\|\big]^T$ and invoking the results of Lemmas \ref{Lemma_x}, \ref{Lemma_r}, and \ref{Lemma_s}, we have the following inequality for $k \geq K+\bar{N}-1$
\begin{equation}
	\begin{aligned}
		\begin{bmatrix}
			\boldsymbol{\xi}^{k+1}\\
			\vdots\\
			\boldsymbol{\xi}^{k-\bar{N}+2}
		\end{bmatrix}
	 \leq \underbrace{(\mathbf{M}^1+\lambda\mathbf{M}^2)}_ {\mathbf{M}(\lambda)}
	 \begin{bmatrix}
		 \boldsymbol{\xi}^{k}\\
		 \vdots\\
		 \boldsymbol{\xi}^{k-\bar{N}+1}
	 \end{bmatrix}	
	\end{aligned}
\end{equation}
where $\mathbf{M}^1$ and $\mathbf{M}^2$ are given by
\begin{equation}
	\begin{aligned}
		& \ \, \mathbf{M}^1 = \qquad \qquad \qquad \qquad \quad \ \mathbf{M}^2 = \\
		&
		\begin{bmatrix}
			\mathbf{M}^1_{a} & \mathbf{M}^1_{b} & \cdots & \mathbf{M}^1_{b} & \mathbf{M}^1_{c}\\
			\mathbf{I}_3 & & & & \\
			& & \ddots & & \\
			& & & \mathbf{I}_3 & \\
		\end{bmatrix},
	 \quad \
		\begin{bmatrix}
			\mathbf{M}^2_{a} & \mathbf{M}^2_{b} & \cdots & \mathbf{M}^2_{b} & \mathbf{M}^2_{c}\\
			\mathbf{0}_{3 \times 3} & & & & \\
			& & \ddots & & \\
			& & & \mathbf{0}_{3 \times 3} & \\
		\end{bmatrix}
	\end{aligned}
\end{equation}
respectively, with
\begin{equation}
\begin{aligned}
\mathbf{M}^1_{a} & =
\begin{bmatrix}
0 & 0 & 0\\
0 & 1 & 0\\
2t & 0 & 0\\
\end{bmatrix}, 
\ \ \quad \mathbf{M}^2_{a} =
\begin{bmatrix}
	n \bar{\beta} q & n \bar{\beta} q & q\\
	n \bar{\beta} & -\eta^{n-1}m & n\\
	n\bar{\beta} t & n\bar{\beta} t & t\\
\end{bmatrix} \\
\mathbf{M}^1_{b} & =
\begin{bmatrix}
0 & 0 & 0\\
0 & 0 & 0\\
2t & 0 & 0
\end{bmatrix},
\ \ \quad \mathbf{M}^2_{b} =
\begin{bmatrix}
	n \bar{\beta} q & n \bar{\beta} q & q\\
	0 & 0 & 0\\
	n\bar{\beta} t & n\bar{\beta} t & t\\
\end{bmatrix} \\
\mathbf{M}^1_{c} & =
\begin{bmatrix}
r_R & 0 & 0\\
0 & 0 & 0\\
2t & 0 & r_P
\end{bmatrix}, 
\quad \mathbf{M}^2_{c} =
\begin{bmatrix}
	n \bar{\beta} q & n \bar{\beta} q & q\\
	0 & 0 & 0\\
	n\bar{\beta} t & n\bar{\beta} t & t\\
\end{bmatrix}\\
t & = {n^2\bar{\beta}Q_P}/{\eta^{n-1}}, \quad q=Q_R\sqrt{n}, \quad m={\alpha_F}/{n}
\end{aligned}
\end{equation}

To establish the $R$-linear convergence of our algorithm, it is sufficient to show that the spectral radius of $\mathbf{M}(\lambda)$ is less than $1$ when $\lambda$ is small enough. To this end, we first prove that the spectral radius of $\mathbf{M}^1$ is $1$, i.e., $\rho(\mathbf{M}^1)=1$. Following Theorem 3.2 in \cite{powell2011calculating}, the determinant of $z\mathbf{I}-\mathbf{M}^1$ is given by $\det(z\mathbf{I}-\mathbf{M}^1)=(z^{\bar{N}}-r_R)(z^{\bar{N}}-r_P)(z-1)z^{\bar{N}-1}$. Since $r_R, r_P \in (0, \, 1)$ holds, we have $\rho(\mathbf{M}^1)=1$. Moreover, the eigenvalue $1$ is simple, and its corresponding right and left eigenvectors are $\mathbf{u} =\mathbf{1}_{\bar{N}} \otimes [0 \ 1 \ 0]^T$ and $\mathbf{w} =[0 \ 1 \ 0 \ \cdots \ 0]^T$, respectively. Denote the simple eigenvalue of $\mathbf{M}(\lambda)$ as a function of $\lambda$, i.e., $p(\lambda)$. Given $\mathbf{M}(\lambda)=\mathbf{M}^1+\lambda \mathbf{M}^2$, we have $p(0)=1$. Using the matrix perturbation theory in Theorem 6.3.12 in \cite{horn2012matrix}, one can obtain $\frac{dp(\lambda)}{d\lambda}\Big|_{\lambda=0}=\frac{\mathbf{w}^T \mathbf{M}^2 \mathbf{u}}{\mathbf{w}^T\mathbf{u}} = -n^{-1}\eta^{n-1} \alpha_F<0$. Since eigenvalues are continuous functions of the elements of a matrix, the simple eigenvalue $p(\lambda)$ is strictly less than $1$ when $\lambda$ is sufficiently small, implying that the spectral radius of $\mathbf{M}(\lambda)$ is less than $1$ when $\lambda$ is sufficiently small. Noting that $\mathbf{M}(\lambda)$ has nonnegative entries and $\mathbf{M}(\lambda)^{\bar{N}+1}$ has all positive entries, from Theorem 8.5.1 and Theorem 8.5.2 in \cite{horn2012matrix}, we can obtain that each entry of $\mathbf{M}(\lambda)^{k}$ will converge to zero at the rate of $\mathcal{O}(\rho(\mathbf{M}(\lambda))^k)$. Therefore, when $\lambda$ is sufficiently small, $\|\mathbf{x}^k - \mathbf{1}_n \otimes \mathbf{x}^*\|$ converges to $0$ at the rate of $\mathcal{O}(\rho(\mathbf{M}(\lambda))^k)$, implying that our algorithm has $R$-linear convergence. \hfill{$\blacksquare$}

\begin{Remark}\label{K_convergence}
	Although adding randomness in optimization parameters in the first $K$ iterations may delay the convergence of our algorithm (as shown in our numerical simulations in Fig. \ref{fig_comparison_convergence_other_algorithms}), it has no influence on the $R$-linear convergence rate. The added randomness in optimization parameters in the first $K$ iterations is key to enable privacy protection, as elaborated in the following subsection.
\end{Remark}

\subsection{Privacy Analysis}

To protect $f_i(\cdot)$, it suffices to protect the gradient function $\nabla f_i(\cdot)$. In decentralized optimization, every participating node receives messages from its neighbors, and hence could exploit them to infer other nodes' private gradients. We denote all information accessible to node $j$ as $\mathcal{I}_j$, which contains the coupling weights, stepsizes, and states associated with node $j$, sent information from node $j$ to its out-neighbors, and received information from node $j$'s in-neighbors to itself. A mathematical representation of $\mathcal{I}_j$ is given as follows:
\begin{equation}\label{collected_information}
	\begin{aligned}
		\mathcal{I}_j = & \big\{ \mathcal{I}_j^{\text{state}}(k) \cup \mathcal{I}_j^{\text{send}}(k) \cup \mathcal{I}_j^{\text{receive}}(k) \, \big| \, k=0,1,\ldots \big\}\\
		& \cup \big\{ \boldsymbol{\Lambda}_j^{k}, \mathbf{R}_{jj}^{k}, \mathbf{A}_{jj}^{k}, \mathbf{C}_{jj}^{k}, \mathbf{B}_{jj}^{k} \big| \, k <K \big\} \\
		& \cup \big\{ \mathbf{R}_{jl}^{k}, \mathbf{A}_{jl}^{k} \big| \, k <K, \forall \, l \in \mathcal{N}_j^{in} \big\} \\
		& \cup \big\{ \mathbf{C}_{mj}^{k}, \mathbf{B}_{mj}^{k} \big| \, k <K, \forall \, m \in \mathcal{N}_j^{out} \big\} \\
		& \cup \big\{ \boldsymbol{\Lambda}_l^{k}, \mathbf{R}_{lm}^{k}, \mathbf{A}_{lm}^{k}, \mathbf{C}_{lm}^{k}, \mathbf{B}_{lm}^{k} \big| \, k \geq K, \forall \, l,m \in \mathcal{V} \big\} \\
	\end{aligned}
\end{equation}
where
\begin{equation}\label{collected_information_detail}
	\begin{aligned}
		& \mathcal{I}_j^{\text{state}}(k) = \big\{\mathbf{x}_j^{k}, \mathbf{y}_j^{k}, \boldsymbol{\Lambda}_j^{k}\mathbf{y}_j^{k}, \mathbf{C}_{jj}^{k} \mathbf{y}_j^{k} + \mathbf{B}_{jj}^{k} ( \nabla f_j^{k+1} - \nabla f_j^{k}) \big\}\\
		& \mathcal{I}_j^{\text{send}}(k) = \big\{ \mathbf{x}_j^{k}, \boldsymbol{\Lambda}_j^{k}\mathbf{y}_j^{k}, \\
		& \qquad \qquad \quad \mathbf{C}_{mj}^{k} \mathbf{y}_j^{k} + \mathbf{B}_{mj}^{k} ( \nabla f_j^{k+1} - \nabla f_j^{k}) \big| \, \forall \, m \in \mathcal{N}_j^{out} \big\}\\
		& \mathcal{I}_j^{\text{receive}}(k) = \big\{ \mathbf{x}_l^{k}, \boldsymbol{\Lambda}_l^{k}\mathbf{y}_l^{k}, \\
		& \qquad \qquad \qquad \mathbf{C}_{jl}^{k} \mathbf{y}_l^{k} + \mathbf{B}_{jl}^{k} ( \nabla f_l^{k+1} - \nabla f_l^{k}) \big| \, \forall \, l \in \mathcal{N}_j^{in} \big\}\\
	\end{aligned}
\end{equation}
represent the respective state information, sent information, and received information of node $j$ at iteration $k$.

Using the defined information set, we are in position to introduce the attacker model and our definition of privacy.

\begin{Definition}\label{honest-but-curious-definition}
	A node $j$ is called honest-but-curious if it follows all protocol steps correctly but is curious and tries to infer the gradient functions of other participating nodes using the information $\mathcal{I}_j$ accessible to itself.
\end{Definition}

\begin{Definition}\label{privacy-preservation-definition}
	For a decentralized network of $n$ nodes, the privacy of node $i$ is preserved if honest-but-curious adversaries cannot distinguish the actual gradient $\nabla f_i(\cdot)$ of node $i$ from an arbitrarily large variation of the gradient $\nabla \tilde{f}_i(\cdot)=\nabla f_i(\cdot) + \boldsymbol{\delta}$ where $\boldsymbol{\delta}$ can be any vector in $\mathbb{R}^d$. That is to say, there exist feasible coupling weights and step-sizes satisfying the requirements in Table II that make the information accessible to honest-but-curious adversaries exactly unchanged under arbitrary variations of node $i$'s gradient, which results in the unidentifiability of node $i$'s actual gradient.
\end{Definition}

Definition \ref{privacy-preservation-definition} means that an honest-but-curious adversary cannot even identify a range for a private function's values and thus is more stringent than the privacy definition used in unobservability based approaches in \cite{observability1, observability2} which define privacy as the inability of an adversary to uniquely determine a protected value. Note that opacity has also been used to protect information in discrete-event systems \cite{saboori2013verification, lefebvre2020privacy}, which was recently extended to continuous-state dynamical systems \cite{ramasubramanian2019notions}. In contrast to unobservability based approaches which consider adversaries having access to the entire output and control histories, opacity based approaches only consider adversaries having access to snapshots of the output and the set of controls. So opacity based privacy is weaker than unobservability based privacy \cite{ramasubramanian2019notions}, and hence is also weaker than the privacy defined here. Furthermore, in contrast to unobservability and opacity based approaches focusing on protecting some state values, we protect participating nodes' gradient function values over the entire domain (or protecting both function types and function parameters), which is more challenging.

Before analyzing the privacy-preserving performance of our proposed algorithm, we first show that using deterministic parameters may easily cause privacy breaches in existing decentralized algorithms. We take the $\mathcal{A}\mathcal{B}$ algorithm in \cite{xin2018linear} as an example. In the $\mathcal{A}\mathcal{B}$ algorithm, node $i$ updates its states $\mathbf{x}_i^{k}$ and $\mathbf{y}_i^{k}$ as follows:
\begin{equation}
	\begin{aligned}
		\mathbf{x}_i^{k+1} & = {\textstyle \sum_{j\in \mathcal{N}_i^{in} \cup \{i\}} } r_{ij} \mathbf{x}_j^{k} - \lambda \mathbf{y}_i^{k} \\
		\mathbf{y}_i^{k+1} & = {\textstyle \sum_{j\in \mathcal{N}_i^{in} \cup \{i\}}} c_{ij} \big( \mathbf{y}_j^{k} + \nabla f_j^{k+1} - \nabla f_j^{k} \big)
	\end{aligned}
\end{equation}
At the initial iteration $k=0$, node $i$ computes $\mathbf{x}_i^{1}$ and $\nabla f_i^{1}$, and then sends $c_{ji} \big( \mathbf{y}_i^{0} + \nabla f_i^{1} - \nabla f_i^{0} \big) = c_{ji}\nabla f_i^{1}$ to its out-neighbor $j$ where $\mathbf{y}_i^{0}=\nabla f_i^{0}$. At iteration $k=1$, node $i$ sends $\mathbf{x}_i^{1}$ to its out-neighbor $j$. Note that in \cite{xin2018linear}, node $i$ sets its corresponding coupling weights $c_{ji}$ as $c_{ji}=1/(|\mathcal{N}_i^{out}|+1)$ where $|\mathcal{N}_i^{out}|$ represents the number of node $i$'s out-neighbors. As a result, using $c_{ji}\nabla f_i^{1}$ obtained at $k=0$ and state $\mathbf{x}_i^1$ received at $k=1$, node $j$ can uniquely determine the gradient of node $i$ at $\mathbf{x}_i^{1}$ as long as it knows the number of node $i$'s out-neighbors. Therefore, only employing deterministic or constant parameters makes the gradients of participating nodes easily inferable by their respective neighboring nodes. Following a similar argument, one can also see that other commonly used algorithms have the same issue of breaching individual nodes' privacy.

Next, we show that by adding randomness in interaction parameters, our algorithm can protect the privacy of gradients. We consider two different scenarios: 1) a single adversary acting on its own (i.e., without collusion with other adversaries) and 2) multiple adversaries colluding with each other.

\subsubsection{Non-colluding case}

In this case, we assume that there are $1 \leq q \leq n-1$ honest-but-curious nodes which try to infer node $i$'s gradient function without sharing information with each other.

\begin{Theorem}\label{single-case-theorem}
	In our decentralized optimization algorithm, the privacy of node $i$ can be preserved if $|\mathcal{N}_i^{out} \cup \mathcal{N}_i^{in}| \geq 2$ holds.
\end{Theorem}

{\it Proof}: Our idea is based on the indistinguishability of $\nabla f_i(\cdot)$'s {\it arbitrary} variations to any honest-but-curious node $j$. According to Definitions \ref{honest-but-curious-definition} and \ref{privacy-preservation-definition}, all the information available to node $j$ to infer $\nabla f_i(\cdot)$ is $\mathcal{I}_j$ (cf. (\ref{collected_information})), and we have to prove that when $\nabla f_i(\cdot)$ is altered to $\nabla \tilde{f}_i(\cdot)$ (could have an {\it arbitrarily} large difference from $\nabla f_i(\cdot)$), the information accessible to node $j$, i.e., $\tilde{\mathcal{I}}_j$, could be exactly the same as $\mathcal{I}_j$ in (\ref{collected_information}). Therefore, we only need to prove that there exists such $\nabla \tilde{f}_i(\cdot)$ that makes $\tilde{\mathcal{I}}_j = \mathcal{I}_j$ hold. Given $|\mathcal{N}_i^{out} \cup \mathcal{N}_i^{in}| \geq 2$, there must exist a node $m \in \mathcal{N}_i^{out} \cup \mathcal{N}_i^{in}$ such that $m \neq j$ holds. So we only need to show that there exist feasible parameters (coupling weights and stepsizes satisfying the requirements in Table II) making $\tilde{\mathcal{I}}_j = \mathcal{I}_j$ hold under $\nabla \tilde{f}_i(\cdot) = \nabla f_i(\cdot) + \boldsymbol{\delta}$ and $\nabla \tilde{f}_m(\cdot) = \nabla f_m(\cdot) - \boldsymbol{\delta}$ for {\it any} $\boldsymbol{\delta}= [\delta_1, \cdots, \delta_d]^T \in \mathbb{R}^d$.

We consider $m \in \mathcal{N}_i^{out}$ and $m \in \mathcal{N}_i^{in}$, separately (note that if $m \in \mathcal{N}_i^{out} \cap \mathcal{N}_i^{in}$ holds, either of the considered cases can be used in the argument to draw a same conclusion):

Case I: If $m \in \mathcal{N}_i^{out}$ holds, it can be proven that $\tilde{\mathcal{I}}_j = \mathcal{I}_j$ holds for any $\boldsymbol{\delta} \in \mathbb{R}^d$ under $\nabla \tilde{f}_i(\cdot) = \nabla f_i(\cdot) + \boldsymbol{\delta}$, $\nabla \tilde{f}_m(\cdot) = \nabla f_m(\cdot) - \boldsymbol{\delta}$, and the following parameters
\begin{equation}\label{proper-coupling-weights-I}
	\left\lbrace\begin{aligned}
		& \mathbf{\tilde{\Lambda}}_{p}^{0} = \mathbf{\Lambda}_{p}^{0} \quad \forall \, p \in \mathcal{V} \setminus \{i, \, m\} \\
		& \tilde{\lambda}_{i}^{0}(l) = \lambda_{i}^{0}(l) \mathbf{y}_i^{0}[l] / (\mathbf{y}_i^{0}[l]+\delta_l) \\
		& \tilde{\lambda}_{m}^{0}(l) = \lambda_{m}^{0}(l) \mathbf{y}_m^{0}[l] / (\mathbf{y}_m^{0}[l]-\delta_l) \\
		& \tilde{\mathbf{R}}_{pq}^{0} = \mathbf{R}_{pq}^{0} \quad \forall \, p,q \in \mathcal{V}\\
		& \tilde{\mathbf{A}}_{pq}^{0} = \mathbf{A}_{pq}^{0} \quad \forall \, p,q \in \mathcal{V} \\	
		& \tilde{\mathbf{C}}_{pq}^{0} = \mathbf{C}_{pq}^{0} \quad \forall \, p,q \in \mathcal{V} \setminus \{i, \, m\} \\
		& \tilde{c}_{pi}^{0}(l)= c_{pi}^{0}(l) \mathbf{y}_i^{0}[l] / (\mathbf{y}_i^{0}[l]+\delta_l) \quad \forall \, p \in \mathcal{V} \setminus \{m\} \\
		& \tilde{c}_{mi}^{0}(l)=(c_{mi}^{0}(l) \mathbf{y}_i^{0}[l] + \delta_l)/ (\mathbf{y}_i^{0}[l]+\delta_l) \\
		& \tilde{c}_{pm}^{0}(l)= c_{pm}^{0}(l) \mathbf{y}_m^{0}[l] / (\mathbf{y}_m^{0}[l]-\delta_l) \quad \forall \, p \in \mathcal{V} \setminus \{m\} \\
		& \tilde{c}_{mm}^{0}(l)=(c_{mm}^{0}(l) \mathbf{y}_m^{0}[l] - \delta_l)/ (\mathbf{y}_m^{0}[l]-\delta_l) \\
		& \tilde{\mathbf{C}}_{pq}^{0} = \mathbf{C}_{pq}^{0} \quad \forall \, p,q \in \mathcal{V} \\
		& \tilde{\mathbf{B}}_{pq}^{0} = \mathbf{B}_{pq}^{0} \quad \forall \, p,q \in \mathcal{V} \\
		& \mathbf{\tilde{\Lambda}}_{p}^{k} = \mathbf{\Lambda}_{p}^{k} \quad \forall \, p \in \mathcal{V}, \, k=1, 2, \ldots \\
		& \tilde{\mathbf{R}}_{pq}^{k} = \mathbf{R}_{pq}^{k} \quad \forall \, p,q \in \mathcal{V}, \, k=1, 2, \ldots \\
		& \tilde{\mathbf{A}}_{pq}^{k} = \mathbf{A}_{pq}^{k} \quad \forall \, p,q \in \mathcal{V}, \, k=1, 2, \ldots \\
		& \tilde{\mathbf{C}}_{pq}^{k} = \mathbf{C}_{pq}^{k} \quad \forall \, p,q \in \mathcal{V}, \, k=1, 2, \ldots \\
		& \tilde{\mathbf{B}}_{pq}^{k} = \mathbf{B}_{pq}^{k} \quad \forall \, p,q \in \mathcal{V}, \, k=1, 2, \ldots \\
\end{aligned}\right.
\end{equation}
where $l=1,\cdots,d$ and ``$\setminus$'' represents set subtraction.

Case II: If $m \in \mathcal{N}_i^{in}$ holds, it can be verified that $\tilde{\mathcal{I}}_j = \mathcal{I}_j$ is true for any $\boldsymbol{\delta} \in \mathbb{R}^d$ under $\nabla \tilde{f}_i(\cdot) = \nabla f_i(\cdot) + \boldsymbol{\delta}$, $\nabla \tilde{f}_m(\cdot) = \nabla f_m(\cdot) - \boldsymbol{\delta}$, and the following parameters
\begin{equation}\label{proper-coupling-weights-II}
	\left\lbrace\begin{aligned}
		& \mathbf{\tilde{\Lambda}}_{p}^{0} = \mathbf{\Lambda}_{p}^{0} \quad \forall \, p \in \mathcal{V} \setminus \{i, \, m\} \\
		& \tilde{\lambda}_{i}^{0}(l) = \lambda_{i}^{0}(l) \mathbf{y}_i^{0}[l] / (\mathbf{y}_i^{0}[l]+\delta_l) \\
		& \tilde{\lambda}_{m}^{0}(l) = \lambda_{m}^{0}(l) \mathbf{y}_m^{0}[l] / (\mathbf{y}_m^{0}[l]-\delta_l) \\
		& \tilde{\mathbf{R}}_{pq}^{0} = \mathbf{R}_{pq}^{0} \quad \forall \, p,q \in \mathcal{V} \\
		& \tilde{\mathbf{A}}_{pq}^{0} = \mathbf{A}_{pq}^{0} \quad \forall \, p,q \in \mathcal{V} \\
		& \tilde{\mathbf{C}}_{pq}^{0} = \mathbf{C}_{pq}^{0} \quad \forall \, p,q \in \mathcal{V} \setminus \{i, \, m\} \\
		& \tilde{c}_{pi}^{0}(l)= c_{pi}^{0}(l) \mathbf{y}_i^{0}[l] / (\mathbf{y}_i^{0}[l]+\delta_l) \quad \forall \, p \in \mathcal{V} \setminus \{m\} \\
		& \tilde{c}_{ii}^{0}(l)=(c_{ii}^{0}(l) \mathbf{y}_i^{0}[l] + \delta_l)/ (\mathbf{y}_i^{0}[l]+\delta_l) \\
		& \tilde{c}_{pm}^{0}(l)= c_{pm}^{0}(l) \mathbf{y}_m^{0}[l] / (\mathbf{y}_m^{0}[l]-\delta_l) \quad \forall \, p \in \mathcal{V} \setminus \{m\} \\
		& \tilde{c}_{im}^{0}(l)=(c_{im}^{0}(l) \mathbf{y}_m^{0}[l] - \delta_l)/ (\mathbf{y}_m^{0}[l]-\delta_l) \\
		& \tilde{\mathbf{C}}_{pq}^{0} = \mathbf{C}_{pq}^{0} \quad \forall \, p,q \in \mathcal{V} \\
		& \tilde{\mathbf{B}}_{pq}^{0} = \mathbf{B}_{pq}^{0} \quad \forall \, p,q \in \mathcal{V} \\
		& \mathbf{\tilde{\Lambda}}_{p}^{k} = \mathbf{\Lambda}_{p}^{k} \quad \forall \, p \in \mathcal{V}, \, k=1, 2, \ldots \\
		& \tilde{\mathbf{R}}_{pq}^{k} = \mathbf{R}_{pq}^{k} \quad \forall \, p,q \in \mathcal{V}, \, k=1, 2, \ldots \\
		& \tilde{\mathbf{A}}_{pq}^{k} = \mathbf{A}_{pq}^{k} \quad \forall \, p,q \in \mathcal{V}, \, k=1, 2, \ldots \\
		& \tilde{\mathbf{C}}_{pq}^{k} = \mathbf{C}_{pq}^{k} \quad \forall \, p,q \in \mathcal{V}, \, k=1, 2, \ldots \\
		& \tilde{\mathbf{B}}_{pq}^{k} = \mathbf{B}_{pq}^{k} \quad \forall \, p,q \in \mathcal{V}, \, k=1, 2, \ldots \\
	\end{aligned}\right.
\end{equation}
where $l=1,\cdots,d$.

To see why the parameter setting in (\ref{proper-coupling-weights-I}) can ensure $\tilde{\mathcal{I}}_j = \mathcal{I}_j$ in the case $m\in \mathcal{N}_i^{out}$, we consider $k=0$ and $k\geq 1$, separately. Under the feasible parameters for $k=0$, it can be verified that the information accessible to node $j$ at the initial iteration $k=0$ keeps unchanged and the states of each node $p \in \mathcal{V}$ satisfy $\mathbf{\tilde{x}}_p^1 = \mathbf{x}_p^1$ and $\mathbf{\tilde{y}}_p^1 = \mathbf{y}_p^1$. Then by setting feasible parameters for iterations $k\geq 1$ the same as the original ones without gradient variations, it is apparent that the information accessible to node $j$ keeps unchanged for iterations $k \geq 1$. Following a similar argument, one can verify that the parameter setting in (\ref{proper-coupling-weights-II}) can ensure $\tilde{\mathcal{I}}_j = \mathcal{I}_j$ in the case $m\in \mathcal{N}_i^{in}$.

In summary, we have $\tilde{\mathcal{I}}_j = \mathcal{I}_j$ for $\nabla \tilde{f}_i(\cdot) = \nabla f_i(\cdot) + \boldsymbol{\delta}$ under {\it any} $\boldsymbol{\delta} \in \mathbb{R}^d$, which means that node $j$ cannot distinguish arbitrarily large variations on $\nabla f_i(\cdot)$ based on its accessible information. Therefore, our proposed algorithm can protect the privacy of node $i$ if $|\mathcal{N}_i^{out} \cup \mathcal{N}_i^{in}| \geq 2$ holds. \hfill{$\blacksquare$}

\begin{Remark}
	Note that if we view the feasible parameters as solutions to guaranteeing $\tilde{\mathcal{I}}_j = \mathcal{I}_j$, then there exist infinitely many solutions. The proof of Theorem \ref{single-case-theorem} just provides one such solution.
\end{Remark}

Next we show that if the condition in Theorem \ref{single-case-theorem} is not met, then the privacy of node $i$ can be breached. 

\begin{Theorem}\label{single-case-theorem-2}
	In our decentralized optimization algorithm, the privacy of node $i$ cannot be preserved against node $j$ if node $j$ is the only in-neighbor and out-neighbor of node $i$, i.e., $\mathcal{N}_i^{out} = \mathcal{N}_i^{in}= \{j\}$. 
\end{Theorem}

{\it Proof}: When $\mathcal{N}_i^{out} = \mathcal{N}_i^{in}= \{j\}$ holds, one can get the dynamics of $\mathbf{y}_i^{k}$ from (\ref{update-rule-02}) as follows:
\begin{equation}\label{y-dynamics-1}
\begin{aligned}
\mathbf{y}_i^{k+1} & = \mathbf{C}_{ii}^{k} \mathbf{y}_i^{k} + \mathbf{B}_{ii}^{k} \big( \nabla f_i^{k+1} - \nabla f_i^{k} \big) + \mathbf{C}_{ij}^{k} \mathbf{y}_j^{k} \\
& \qquad + \mathbf{B}_{ij}^{k} \big( \nabla f_j^{k+1} - \nabla f_j^{k} \big)
\end{aligned}
\end{equation}	
According to the parameter design in Table II, we have $\mathbf{y}_i^{k} = \mathbf{C}_{ii}^{k} \mathbf{y}_i^{k} + \mathbf{C}_{ji}^{k} \mathbf{y}_i^{k}$ and $\nabla f_i^{k+1} - \nabla f_i^{k} = (\mathbf{B}_{ii}^{k}+ \mathbf{B}_{ji}^{k}) (\nabla f_i^{k+1} - \nabla f_i^{k})$ for $k \in \mathbb{Z}_{\geq 0}$ based on the facts $\mathbf{C}_{ii}^{k} + \mathbf{C}_{ji}^{k} = \mathbf{I}_d$ and $\mathbf{B}_{ii}^{k} + \mathbf{B}_{ji}^{k} = \mathbf{I}_d$. So we can rewrite (\ref{y-dynamics-1}) as
\begin{equation}\label{y-dynamics}
	\begin{aligned}
	\mathbf{y}_i^{k+1} = & \mathbf{y}_i^{k} - \mathbf{C}_{ji}^{k} \mathbf{y}_i^{k} + (\nabla f_i^{k+1} - \nabla f_i^{k}) + \mathbf{C}_{ij}^{k} \mathbf{y}_j^{k} \\
	& - \mathbf{B}_{ji}^{k} (\nabla f_i^{k+1} - \nabla f_i^{k}) + \mathbf{B}_{ij}^{k} (\nabla f_j^{k+1} - \nabla f_j^{k})
	\end{aligned}
\end{equation}
Denote $\mathbf{m}_j^k$ as 
\begin{equation}\label{privacy_proof_1}
\begin{aligned}
\mathbf{m}_j^k & = - \mathbf{C}_{ji}^{k} \mathbf{y}_i^{k} + \mathbf{C}_{ij}^{k} \mathbf{y}_j^{k} - \mathbf{B}_{ji}^{k} (\nabla f_i^{k+1} - \nabla f_i^{k})\\
& \qquad + \mathbf{B}_{ij}^{k} (\nabla f_j^{k+1} - \nabla f_j^{k})
\end{aligned}
\end{equation}
Note that $\mathbf{m}_j^k$ is accessible to the honest-but-curious node $j$ because $\mathbf{C}_{ji}^{k} \mathbf{y}_i^{k}+\mathbf{B}_{ji}^{k} (\nabla f_i^{k+1} - \nabla f_i^{k})$ is the information node $i$ sends to node $j$, and $\mathbf{C}_{ij}^{k} \mathbf{y}_j^{k} + \mathbf{B}_{ij}^{k} (\nabla f_j^{k+1} - \nabla f_j^{k})$ is the information computed by node $j$. Plugging (\ref{privacy_proof_1}) into (\ref{y-dynamics}), node $j$ can obtain $\mathbf{y}_i^{k+1} - \mathbf{y}_i^{k} = \nabla f_i^{k+1} - \nabla f_i^{k} + \mathbf{m}_j^k$ and further (note $\mathbf{y}_i^{0} = \nabla f_i^{0}$)
\begin{equation}\label{y-dynamics_3}
\begin{aligned}
\mathbf{y}_i^{k+1} = \nabla f_i^{k+1} + {\textstyle \sum_{l=0}^{k}}\mathbf{m}_j^l
\end{aligned}
\end{equation}
Since $\lim_{k \rightarrow \infty}\mathbf{y}_i^{k+1} = \mathbf{0}_d$ holds as $k$ goes to infinity, we have $\lim_{k \rightarrow \infty} \nabla f_i^{k+1} = \lim_{k \rightarrow \infty} \nabla f_i(\mathbf{x}_i^{k+1}) = - \lim_{k \rightarrow \infty} \sum_{l=0}^{k}\mathbf{m}_j^l$. Given $\lim_{k \rightarrow \infty} \mathbf{x}_i^{k+1} = \lim_{k \rightarrow \infty} \mathbf{x}_j^{k+1}$ $= \mathbf{x}^{*}$, node $j$ also knows $\lim_{k \rightarrow \infty} \mathbf{x}_i^{k+1}$, meaning that an honest-but-curious node $j$ can infer the gradient of node $i$ at the global optimal solution $\mathbf{x}^{*}$ using (\ref{y-dynamics_3}). Therefore, the privacy of node $i$ cannot be preserved against node $j$ when node $j$ is the only in-neighbor and out-neighbor of node $i$. \hfill{$\blacksquare$}

\subsubsection{Colluding case}

When a set of honest-but-curious nodes $\mathcal{A}$ collude, they can share information with each other to infer node $i$'s gradient function.

\begin{Theorem}\label{multiple-case-theorem}
	In our decentralized optimization algorithm, the privacy of node $i$ can be preserved against a set of honest-but-curious nodes $\mathcal{A}$ if $(\mathcal{N}_i^{out} \cup \mathcal{N}_i^{in}) \not\subset \mathcal{A}$ holds, i.e., there exits at least one node that belongs to $\mathcal{N}_i^{out} \cup \mathcal{N}_i^{in}$ but not $\mathcal{A}$.
\end{Theorem}

{\it Proof}: To show that the privacy of node $i$ can be preserved, we have to show that no honest-but-curious node $j\in \mathcal{A}$ can distinguish arbitrary variations on $\nabla f_i(\cdot)$. In the colluding case, each node in $\mathcal{A}$ has access to the information accessible to any node in $\mathcal{A}$. So we represent the set of accessible information as
\begin{equation}\label{collected_information_A}
\begin{aligned}
\mathcal{I_A}= \big\{ \mathcal{I}_j \, \big| \, \forall \, j\in \mathcal{A} \big\}\\
\end{aligned}
\end{equation}
where $\mathcal{I}_j$ is given by (\ref{collected_information}). Following the same line of reasoning in Theorem \ref{single-case-theorem}, to prove that node $j\in \mathcal{A}$ cannot distinguish arbitrary variations on $\nabla f_i(\cdot)$, it suffices to prove when $\nabla f_i(\cdot)$ is altered to $\nabla \tilde{f}_i(\cdot)$ with any $\boldsymbol{\delta} \in \mathbb{R}^d$, the information accessible to the set of nodes in $\mathcal{A}$, i.e., $\tilde{\mathcal{I}}_\mathcal{A}$, could be exactly the same as $\mathcal{I_A}$ in (\ref{collected_information_A}).

If $(\mathcal{N}_i^{out} \cup \mathcal{N}_i^{in}) \not\subset \mathcal{A}$ is true, there must exist a node $m \in \mathcal{N}_i^{out} \cup \mathcal{N}_i^{in}$ such that $m \notin \mathcal{A}$ holds. It can be verified that under $\nabla \tilde{f}_i(\cdot) = \nabla f_i(\cdot) + \boldsymbol{\delta}$ and $\nabla \tilde{f}_m(\cdot) = \nabla f_m(\cdot) - \boldsymbol{\delta}$, there exist respective feasible parameters in (\ref{proper-coupling-weights-I}) (for the case $m \in \mathcal{N}_i^{out}$) and (\ref{proper-coupling-weights-II}) (for the case $m \in \mathcal{N}_i^{in}$) that satisfy the requirements in Table II and make $\tilde{\mathcal{I}}_\mathcal{A}=\mathcal{I_A}$ hold under any $\boldsymbol{\delta} \in \mathbb{R}^d$. Therefore, our proposed algorithm can preserve the privacy of node $i$ against the set of honest-but-curious nodes $\mathcal{A}$ if $(\mathcal{N}_i^{out} \cup \mathcal{N}_i^{in}) \not\subset \mathcal{A}$ holds. \hfill{$\blacksquare$}

Next we show that if the condition in Theorem \ref{multiple-case-theorem} is not met, then the privacy of node $i$ can be breached. 

\begin{Theorem}\label{multiple-case-theorem-2}
	In our decentralized optimization algorithm, the privacy of node $i$ cannot be preserved when all the in-neighbors and out-neighbors belong to $\mathcal{A}$, i.e., $(\mathcal{N}_i^{out} \cup \mathcal{N}_i^{in}) \subset \mathcal{A}$.
\end{Theorem}

{\it Proof}: Following a similar line of reasoning for Theorem \ref{single-case-theorem-2}, Theorem \ref{multiple-case-theorem-2} can be easily obtained and hence we omit the proof here. \hfill{$\blacksquare$}

\begin{Remark}
	It is worth noting that even using time-varying parameters, the $\mathcal{AB}$ algorithm cannot guarantee the privacy defined in Definition \ref{privacy-preservation-definition}. This is because the $\mathcal{AB}$ algorithm does not allow negative coupling weights, leading to no feasible solutions for the parameters required in (\ref{proper-coupling-weights-I}) and (\ref{proper-coupling-weights-II}), which are necessary to make adversaries' accessible information unchanged under arbitrary variations on the gradient.
\end{Remark}

\begin{Remark}\label{K_privacy}
	From Theorem 1, one can see that the parameter $K$ does not affect optimization accuracy. Furthermore, from (\ref{proper-coupling-weights-I}) and (\ref{proper-coupling-weights-II}), we can see that only changing the coupling weights/stepsizes in the initial iteration $k = 0$ is enough to cover gradient variations. In other words, any $K\geq 1$ is sufficient to protect the defined privacy for gradients. However, a larger $K$ does provide additional privacy protection to more intermediate states, if the privacy of intermediate states is of concern. This is because each node $i$ sends $\mathbf{\Lambda}_i^{k} \mathbf{y}_i^{k}$ and $\mathbf{C}_{ji}^{k} \mathbf{y}_i^{k} + \mathbf{B}_{ji}^{k} \big( \nabla f_i(\mathbf{x}_i^{k+1}) - \nabla f_i(\mathbf{x}_i^{k})\big)$ to its out-neighbors $j \in \mathcal{N}_i^{out}$. In the first $K$ iterations, as randomness is injected into the coupling weights ($\mathbf{C}_{ji}^{k}$ and $\mathbf{B}_{ji}^{k}$) and stepsizes ($\mathbf{\Lambda}_i^{k}$), the intermediate states $\mathbf{y}_i^{k}$ and $\nabla f_i(\mathbf{x}_i^{k+1}) - \nabla f_i(\mathbf{x}_i^{k})$ are uninferable by receiving nodes (note that in our algorithm $\mathbf{x}_i^0$ is randomly initialized in $\mathbb{R}^d$ and hence does not carry sensitive information). Therefore, if intermediate states in the first $M$ iterations need to be protected, we can set $K$ to $M$ to achieve this goal without compromising accuracy (note that the gradient will always be protected for any $K\geq 1$). It is also worth noting that although the randomness added in the first $K$ iterations does not affect the convergence rate (as proven in Theorem \ref{theorem_convergence}), it does delay the convergence since the algorithm only starts to converge after iteration $K$ (see Fig. \ref{fig_comparison_convergence_other_algorithms}). Of course, if the privacy of intermediate states is not of concern and only gradients need to be protected, it is preferable to set $K=1$ to minimize the delay in the convergence process.
\end{Remark}

\begin{Remark}\label{remark_side_information}
	In the case where an adversary has side information such as a range of coefficients of the objective function, we cannot prove the adversary's inability to distinguish $\nabla f_i(x_i)$ from $\nabla \tilde{f}_i(x_i) = \nabla f_i(x_i) + \boldsymbol{\delta}$ under an arbitrary $\boldsymbol{\delta}$. Although it is difficult to quantify the influence of side information (e.g., a certain range of parameters) on enabled privacy due to possible complicated ways that the side information can affect the gradient, we can show that in specific cases, such as in the rendezvous problem, the adversary's side information on parameter range could not be tightened after running the algorithm. For the sake of simple exposition, we consider the one-dimensional rendezvous problem. Each node $i$ has a local objective function $f_i(x_i)=\frac{1}{2}(x_i-p_i)^2$. The gradient function of node $i$ is $\nabla f_i (x_i) = x_i-p_i$. Following our analysis in Theorem \ref{single-case-theorem} and Theorem \ref{multiple-case-theorem}, the privacy of node $i$ can be preserved if it has an in-neighbor or out-neighbor not colluding with the adversarial node, i.e., the adversary cannot distinguish whether the gradient function of node $i$ is $\nabla f_i (x_i)=x_i-p_i$ or $\nabla \tilde{f}_i(x_i) \triangleq x_i-p_i + \delta = x_i-\tilde{p}_i$ ($\tilde{p}_i \triangleq p_i-\delta$) under any $\delta$ in $\mathbb{R}$. If the adversary has side information of the function type (i.e., quadratic) and parameter range $p_i\in [a,b]$, then the range of $\delta$ has to be reduced to $[p_i-b, p_i-a]$. This is because otherwise $\tilde{p}_i$ in $\nabla \tilde{f}_i(x_i)$ will not be within the initial range $[a,b]$ known to the adversary. However, for any $\delta \in [p_i-b, p_i-a]$ (i.e., $\tilde{p}_i \in [a,b]$), following our analysis in Theorem \ref{single-case-theorem} and Theorem \ref{multiple-case-theorem}, there always exist feasible parameters making all information known to the adversarial node under $\nabla \tilde{f}_i(x_i)$ exactly the same as under $\nabla f_i(x_i)$. Therefore, the adversarial node cannot distinguish $\nabla f_i(x_i)$ from $\nabla \tilde{f}_i(x_i)$ for any $\delta \in [p_i-b, p_i-a]$, meaning that the adversarial node cannot get a tighter range of $p_i$ than the initial knowledge (i.e., $[a,b]$) after running the algorithm.
\end{Remark}

\begin{Remark}\label{remark_eavesdropper}
	Our approach can be extended to enable privacy against adversaries wiretapping all communication links without compromising algorithmic accuracy by patching partially homomorphic encryption. For example, in our prior work \cite{zhang2018admm} and \cite{zhang2018enabling}, we can let each node generate and flood its public key before the optimization iteration starts. Then in decentralized implementation, a node encrypts its messages to be sent, which can be decrypted by a legitimate recipient without the help of any third party.
\end{Remark}

\section{Numerical Simulations}

\subsection{Privacy Protection in the Rendezvous Problem}

First let's consider the distributed rendezvous problem where a group of nodes want to agree on the nearest meeting point without revealing each other's initial position \cite{huang2015differentially}. Mathematically this can be modeled as the problem $\min_{x\in \mathbb{R}^d} \, \ F(\mathbf{x})=\sum_{i=1}^{n} f_i(\mathbf{x}) = \sum_{i=1}^{n} \frac{1}{2} \| \mathbf{x}-\mathbf{p}_i\|^2$, where $\mathbf{p}_i$ represents the initial position of node $i$. 

\begin{figure}
	\begin{center}
		\includegraphics[width=0.33\textwidth]{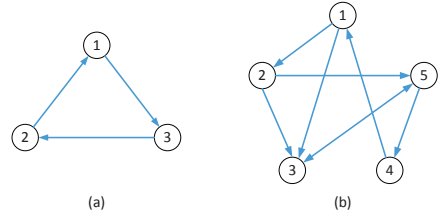}
	\end{center}
	\vspace{-0.3cm}\caption{Graphs used in simulations: (a) a directed cycle graph with $3$ nodes; (b) a strongly connected graph with $5$ nodes.}
	\label{fig_graph_topology}
\end{figure}

For the simplicity of exposition, we consider the $d=1$ case but similar results can be obtained when $d\neq 1$. We consider three nodes connected in a directed cycle as shown in Fig. \ref{fig_graph_topology} (a). Let node $3$ be an honest-but-curious node which collects received data in an attempt to learn the gradient function of node $1$. Node $2$ does not collude with node $3$. Parameter $K$ was set to $3$. In the simulation, we first ran our algorithm and recorded $\mathcal{I}_3$, the information accessible to node $3$ (cf. (\ref{collected_information})). Then we show that information accessible to node $3$ can be exactly the same under a completely different gradient function $\nabla \tilde{f}_1 (\mathbf{x}_1)$.

\begin{figure}
	\begin{center}
		\includegraphics[width=0.45\textwidth]{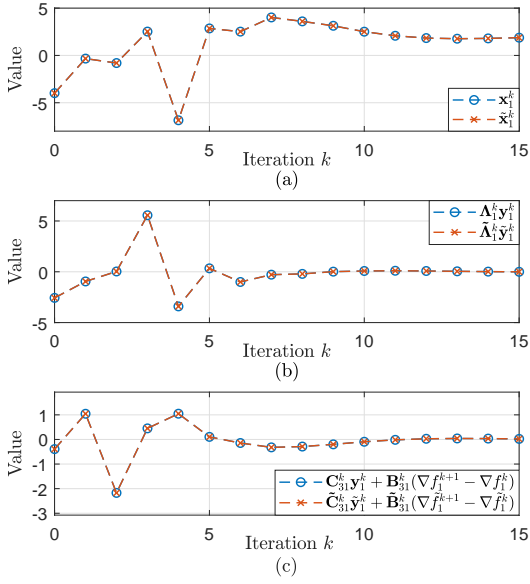}
	\end{center}
	\vspace{-0.4cm}\caption{The information accessible to node $3$ are the same under two different gradient functions of node $1$ depicted in Fig. \ref{fig_gradient_comparison}.}
	\label{fig_exchanged_information_comparison}
\end{figure}

\begin{figure}
	\begin{center}
		\includegraphics[width=0.45\textwidth]{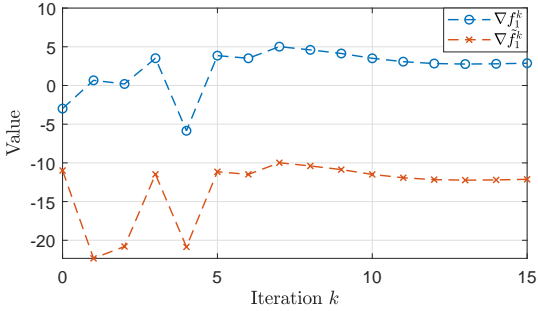}
	\end{center}
	\vspace{-0.2cm}\caption{The two different gradient functions of node $1$ that lead to identical observations at node $3$.}
	\label{fig_gradient_comparison}
\end{figure}

We represent the information accessible to node $3$ in the new implementation as $\tilde{\mathcal{I}}_3$. Fig. \ref{fig_exchanged_information_comparison} shows $\mathbf{x}_{1}^{k}$, $\mathbf{\Lambda}_{1}^{k}\mathbf{y}_{1}^{k}$, and $\mathbf{C}_{31}^{k}\mathbf{y}_{1}^{k}+\mathbf{B}_{31}^{k}( \nabla f_1^{k+1}- \nabla f_1^{k})$ in $\mathcal{I}_3$ and $\tilde{\mathbf{x}}_{1}^{k}$, $\mathbf{\tilde{\Lambda}}_1^k\tilde{\mathbf{y}}_{1}^{k}$, and $\mathbf{\tilde{C}}_{31}^k \tilde{\mathbf{y}}_{1}^{k} + \mathbf{\tilde{B}}_{31}^k ( \nabla \tilde{f}_1^{k+1}- \nabla \tilde{f}_1^{k})$ in $\tilde{\mathcal{I}}_3$, respectively. It can be seen that the trajectory of the observations of node $3$ in both cases are identical. Fig. \ref{fig_gradient_comparison} shows the corresponding $\nabla f_1 (\mathbf{x}_1)$ and $\nabla \tilde{f}_1 (\mathbf{x}_1)$, which are clearly different. Since node $3$ receives the same information under $\nabla \tilde{f}_1 (\mathbf{x}_1) \neq \nabla f_1 (\mathbf{x}_1)$, it has no way to infer the real gradient function of node $1$.

\subsection{Distributed Estimation Problem}

In this subsection we focus on the convergence performance of our algorithm and compare our algorithm with other decentralized optimization algorithms. We consider the canonical distributed estimation problem $\min_{\substack{\mathbf{x} \in \mathbb{R}^d}} \,\ F(\mathbf{x}) = \sum_{i=1}^{n} \big( \| \mathbf{z}_i - \mathbf{Q}_i \mathbf{x} \|^2 + \sigma_i \| \mathbf{x} \|^2 \big)$ in \cite{xu2017convergence}, where $n$ nodes cooperatively measure a certain unknown parameter $\mathbf{x} \in \mathbb{R}^d$. In this problem, each node $i$ has access to its local cost function $f_i(\mathbf{x}_i) = \| \mathbf{z}_i - \mathbf{Q}_i \mathbf{x}_i \|^2 + \sigma_i \| \mathbf{x}_i \|^2$ with $\mathbf{Q}_i \in \mathbb{R}^{s\times d}$ being its measurement matrix and $\mathbf{z}_i \in \mathbb{R}^s$ being its measurement data. The regularization parameter $\sigma_i$ can be set to $0$ (resp. a positive value) to make $f_i$ general convex (resp. strongly convex) under appropriate $\mathbf{Q}_i$.

We considered a network of $n=5$ nodes interacting on a strongly connected graph illustrated in Fig. \ref{fig_graph_topology} (b). Dimension parameters $d$ and $s$ were set to $2$ and $3$, respectively. Parameter $K$ was set to $3$. 

\subsubsection{Comparison with Other Decentralized Optimization Algorithms} 
\begin{figure}
	\begin{center}
		\includegraphics[width=0.45\textwidth]{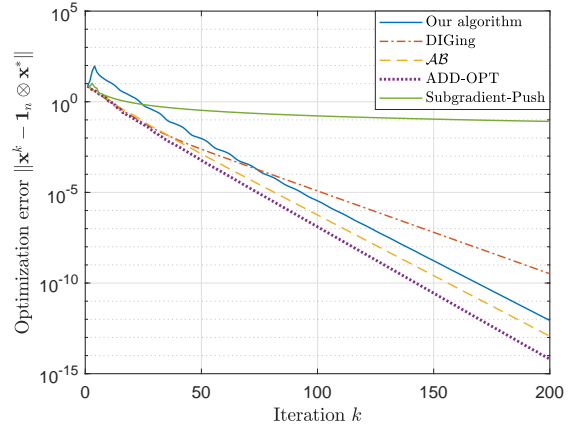}
	\end{center}
	\vspace{-0.2cm}\caption{Comparison of convergence rates of different optimization algorithms.}
	\label{fig_comparison_convergence_other_algorithms}
\end{figure}

We compared our algorithm with DIGing \cite{nedic2017geometrically}, $\mathcal{A}\mathcal{B}$ \cite{xin2018linear}, ADD-OPT \cite{xi2018add}, and Subgradient-Push \cite{nedic2015distributed} to evaluate the influence of privacy design on optimization performance. The stepsize was set to $0.06$ for all of the considered algorithms except Subgradient-Push whose stepsize was set to a diminishing sequence $\lambda^k=1/k$. The simulation results on optimization error $\| \mathbf{x}^k - \mathbf{1}_n \otimes \mathbf{x}^* \|$ are depicted in Fig. \ref{fig_comparison_convergence_other_algorithms}.

From Fig. \ref{fig_comparison_convergence_other_algorithms}, it can be seen that all algorithms with a fixed stepsize had $R$-linear convergence, including our algorithm with guaranteed privacy. This corroborates our statement in Remark \ref{K_convergence} that the added randomness has no influence on the $R$-linear convergence of our algorithm. Of course, it is worth noting that in our algorithm, the privacy-induced randomness delayed convergence to iteration step $k\geq 4$.

\subsubsection{Comparison with \cite{huang2015differentially}} 
\begin{figure}
	\begin{center}
		\includegraphics[width=0.5\textwidth]{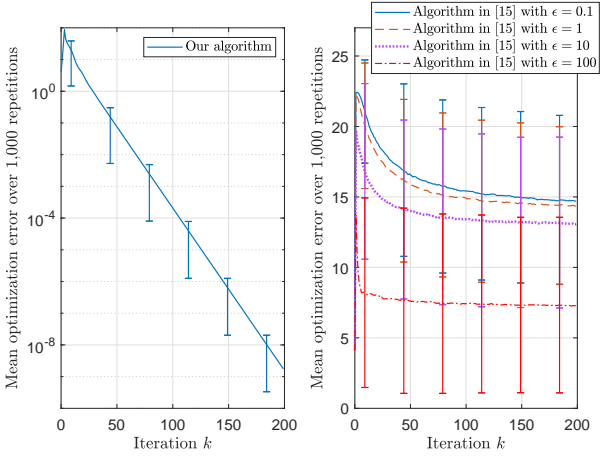}
	\end{center}
	\vspace{-0.2cm}\caption{Comparison of mean optimization error over $1,000$ repetitions between our algorithm and the differential-privacy based approach in \cite{huang2015differentially}.}
	\label{fig_comparison_differential_privacy}
\end{figure}

We also compared our algorithm with the differential privacy based privacy-preserving approach for decentralized optimization proposed in \cite{huang2015differentially}. We ran the algorithm in \cite{huang2015differentially} under four different privacy levels, i.e., $\epsilon = 0.1, 1, 10, 100$, respectively. The domain of optimization in \cite{huang2015differentially} was set to $\mathcal{X}= \{ \mathbf{x} \in \mathbb{R}^2 | \|\mathbf{x}\| \leq 10 \}$. Note that the optimal solution $\mathbf{x}^*=[0.6881, \, 0.5103]^T$ resides in $\mathcal{X}$. For each privacy level $\epsilon$, we repeated the simulation for $1,000$ times, and averaged the optimization error trajectories. Under a stepsize $0.06$, we also measured the mean optimization error of our algorithm over $1,000$ repetitions. The comparison results in Fig. \ref{fig_comparison_differential_privacy} confirm the trade-off between privacy and accuracy for differential-privacy based approaches and demonstrate the advantage of our algorithm in ensuring optimization accuracy.

\subsubsection{Comparison with \cite{lou2018privacy}} 
\begin{figure}
	\begin{center}
		\includegraphics[width=0.45\textwidth]{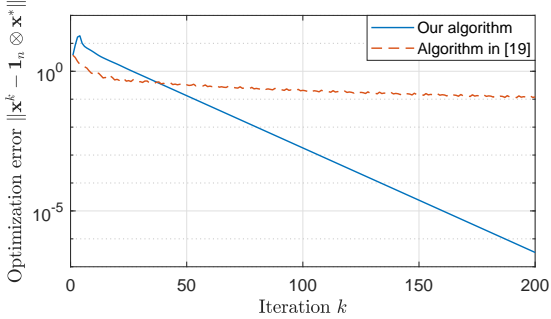}
	\end{center}
	\vspace{-0.2cm}\caption{Comparison of optimization error $\| \mathbf{x}^k - \mathbf{1}_n \otimes \mathbf{x}^* \|$ between our algorithm and the algorithm in \cite{lou2018privacy}.}
	\label{fig_comparison_heter_stepsize_privacy}
\end{figure}

Then we compared our algorithm with the privacy-preserving optimization algorithm in \cite{lou2018privacy}. The closed convex projection set $\mathcal{X}$ was set to $\mathcal{X}= \{ \mathbf{x} \in \mathbb{R}^2 | \|\mathbf{x}\| \leq 10 \}$ for the algorithm in \cite{lou2018privacy}. The other parameters for the algorithm in \cite{lou2018privacy} were set as follows: each node $i$ randomly chose $c_i$ from $(0, \, 5)$, and updated its state with its subgradient once every $T_i$ iterations where $T_i$ was randomly chosen from $\{1,2,3,4,5\}$. The stepsize of our algorithm was set to $0.06$. The results are shown in Fig. \ref{fig_comparison_heter_stepsize_privacy}. Clearly our algorithm has a $R$-linear convergence whereas the convergence rate of \cite{lou2018privacy} is much slower.

\subsubsection{Numerical Simulations on Large-scale Networks} 
\begin{figure}
	\begin{center}
		\includegraphics[width=0.45\textwidth]{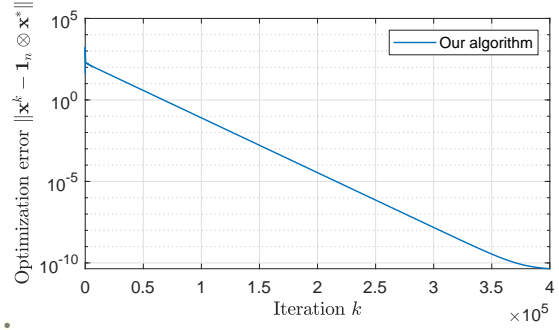}
	\end{center}
	\vspace{-0.2cm}\caption{The optimization error trajectory of our algorithm in a network of $n=100$ nodes.}
	\label{fig_larger_network}
\end{figure}

Finally, we conducted numerical simulations to verify the scalability of our proposed algorithm using a network of $n=100$ nodes. Each agent $i$ was assumed to have two out-neighbors, i.e., $\mathcal{N}_i^{out}=\left\{\overline{i}+1, \, \overline{i+1}+1\right\}$, where the superscript ``$\bar{\quad}$'' represents modulo operation on $n$, i.e., $\overline{i} \triangleq i \mod n$. The evolution of optimization error $\| \mathbf{x}^k - \mathbf{1}_n \otimes \mathbf{x}^* \|$ is shown in Fig. \ref{fig_larger_network}. It can be seen that the convergence rate is still linear, meaning that our proposed algorithm can guarantee the convergence of all nodes to the global optimal solution even when the network size is large.

\section{Conclusions}
In this paper we proposed a dynamics based privacy approach for decentralized optimization. Our approach can enable privacy without compromising optimization accuracy or incurring heavy computation/communication overhead. This is in distinct difference from differential-privacy based approaches which compromise optimization accuracy and encryption based approaches which incur heavy computation/communication overhead. We rigorously characterized the convergence properties of our algorithm and its privacy-preserving performance. In addition, to facilitate the privacy design, we also proposed a general framework of gradient-tracking based decentralization optimization, which includes many commonly used algorithms as special cases. Finally, we provided numerical simulation results to confirm the effectiveness and efficiency of our proposed algorithm.

\section*{Appendix A: Absolute Probability Sequence}

\begin{Definition}\cite{saadatniaki2020decentralized}
	Let $\{\mathbf{A}^k\}$ be a sequence of row-stochastic matrices. A sequence of vectors $\{\boldsymbol{\pi}^k\}$ is an absolute probability sequence for $\{\mathbf{A}^k\}$ if $(\boldsymbol{\pi}^{k})^T=(\boldsymbol{\pi}^{k+1})^T\mathbf{A}^k$ holds for $k\geq 0$.
\end{Definition}

\begin{Definition}\cite{seneta2006non}
	A sequence of row-stochastic matrices $\{\mathbf{A}^k\}$ is ergodic if the following limit exists for $k\geq 0$
	\begin{equation}\label{Def_ergodic}
		\begin{aligned} 
			\lim\limits_{t \rightarrow \infty} \mathbf{A}^t \mathbf{A}^{t-1} \cdots \mathbf{A}^{k+1} \mathbf{A}^{k} =\mathbf{1} (\boldsymbol{\pi}^k)^T
		\end{aligned}
	\end{equation}
	where $\{\boldsymbol{\pi}^k\}$ is an absolute probability sequence for $\{\mathbf{A}^k\}$. 
\end{Definition}

For the absolute probability sequences $\{\mathbf{v}^k\}$ and $\{\boldsymbol{\phi}^k\}$ (associated with $\{\mathbf{\bar{P}}^k\}$ and $\{\mathbf{\bar{R}}^k\}$, respectively) in Sec. 4.2, we have the following lemmas.
\begin{Lemma}\label{Lemma_P}
	Under Assumptions \ref{Assumption_directed_graph} and \ref{Assumption_objective_functions}, and the parameter design in Table II, the sequence of matrices $\{\mathbf{\bar{P}}^k\}$ for $k \geq K$ is row-stochastic and ergodic, and the sequence $\{\mathbf{v}^k\}$ is the unique absolute probability sequence for $\{\mathbf{\bar{P}}^k\}$. Denote $Q_P$ as $Q_P=2n (1+(n\eta^{-n})^{n-1} )/(1-(n^{-1}\eta^n)^{n-1})$ and define a positive integer $N_P$ such that $r_{P} \triangleq Q_P ( 1-(n^{-1}\eta^n)^{n-1} )^{\frac{N_P-1}{n-1}} <1$ holds. Then for any $k\geq K+N_P-1$ and any vector $\mathbf{b}\in \mathbb{R}^{nd}$, vector $\mathbf{a} \triangleq \mathbf{P}^{k} \mathbf{P}^{k-1} \cdots \mathbf{P}^{k-N_P+1}\mathbf{b}$ always satisfies $\|[(\mathbf{I}_n - \mathbf{1}_n (\mathbf{v}^{k+1})^T ) \otimes \mathbf{I}_d] \mathbf{a}\| \leq r_P \|[(\mathbf{I}_n - \mathbf{1}_n (\mathbf{v}^{k-N_P+1})^T) \otimes \mathbf{I}_d ] \mathbf{b}\|$.
\end{Lemma}
{\it Proof}: Given $\mathbf{\bar{P}}^{k} = (\mathbf{\bar{V}}^{k+1})^{-1} \mathbf{\bar{C}}^{k} \mathbf{\bar{V}}^{k}$, it can be verified that $\mathbf{\bar{P}}^{k}$ is row-stochastic since $\sum_{j=1}^{n} p_{ij}^k =1$ holds for all $i=1,\cdots, n$. Using Lemma 4 in \cite{nedic2009distributed}, we can see that the sequence of matrices $\{\mathbf{\bar{P}}^k\}$ for $k \geq K$ is ergodic since the limit in (\ref{Def_ergodic}) exists under Assumptions \ref{Assumption_directed_graph} and \ref{Assumption_objective_functions}, and the parameter design in Table II. Combing $\mathbf{\bar{V}}^k=\text{diag}(\mathbf{v}^k)$ and $\mathbf{\bar{P}}^{k} = (\mathbf{\bar{V}}^{k+1})^{-1} \mathbf{\bar{C}}^{k} \mathbf{\bar{V}}^{k}$ leads to $(\mathbf{v}^k)^T=(\mathbf{v}^{k+1})^T\mathbf{\bar{P}}^{k}$ for $k\geq K$, meaning that $\{\mathbf{v}^k\}$ is the absolute probability sequences for $\{\mathbf{\bar{P}}^k\}$. Further using Lemma 1 in \cite{nedic2016convergence}, we have that $\{\mathbf{v}^k\}$ is the unique absolute probability sequence for $\{\mathbf{\bar{P}}^k\}$ due to the ergodicity of $\{\mathbf{\bar{P}}^k\}$.

Now we show that the positive entries of $\mathbf{\bar{P}}^k$ has a lower bound $\eta^{n}/n$, i.e., $p_{ij}^k \geq \eta^{n}/n$ for any $j \in \mathcal{N}_I^{out} \cup \{i\}$ where $p_{ij}^k$ represents the $ij$-th entry of $\mathbf{\bar{P}}^k$. To this end, we first prove that all entries of $\mathbf{v}^k$ have a lower bound $\eta^{n-1}/n$, i.e., $v_i^k \geq \eta^{n-1}/n$ for $i=1,\cdots,n$ and $k\geq K$. Denoting $\mathbf{\bar{C}}_{k:K}$ as $\mathbf{\bar{C}}_{k:K}= \mathbf{\bar{C}}^k \cdots \mathbf{\bar{C}}^K$, from (\ref{Eqn_vector_v}) we have
\begin{equation}\label{Eqn_v_1}
	\begin{aligned}
		\mathbf{v}^k = \mathbf{\bar{C}}_{k-1:K}\mathbf{v}^K=\mathbf{\bar{C}}_{k-1:K}\mathbf{1}_n/{n}
	\end{aligned}
\end{equation}
for $k\geq K$. Then we divide the proof into two parts: $K \leq k \leq K+n-1$ and $k \geq K+n$, respectively.

\textbf{Part 1:} $v_i^k \geq \eta^{n-1}/n$ for $K \leq k \leq K+n-1$. Given $c_{ii}^k \geq \eta$ for $k\geq K$, one can verify that the inequality $v_i^k \geq \frac{1}{n} [\mathbf{\bar{C}}_{k-1:K}]_{ii} \geq \eta^{k-K}/n \geq {\eta^{n-1}}/{n}$ holds for $i=1,\cdots, n$ and $K \leq k \leq K+n-1$.

\textbf{Part 2:} $v_i^k \geq \eta^{n-1}/n$ for $k \geq K+n$. Under Assumptions \ref{Assumption_directed_graph} and \ref{Assumption_objective_functions}, and the parameter design in Table II, following the arguments in Lemma 2 in \cite{nedic2009distributed}, we can obtain $[\mathbf{\bar{C}_{k-1:k-n+1}}]_{ij} \geq \eta^{n-1}$ for $i, j =1, \cdots, n$. Since $\mathbf{\bar{C}}^k$ is column-stochastic, $\mathbf{\bar{C}}_{k-n:K}$ is also column-stochastic. Further using the fact
$\mathbf{\bar{C}}_{k-1:K} = \mathbf{\bar{C}}_{k-1:K-n+1} \mathbf{\bar{C}}_{k-n:K}$ yields $[\mathbf{\bar{C}}_{k-1:K}]_{ij} \geq \eta^{n-1}$ for $i, j =1, \cdots, n$. Therefore, we have $v_i^k = {1}/{n} \sum_{j=1}^{n} [\mathbf{\bar{C}}_{k-1:K}]_{ij} \geq \eta^{n-1} \geq {\eta^{n-1}}/{n}$ for $i=1,\cdots, n$ and $k \geq K+n$.

In summary, we have $v_i^k \geq \eta^{n-1}/n$ for $i=1,\cdots,n$ and $k\geq K$. Since $\mathbf{\bar{C}}_{k-1:K}$ is column-stochastic, it follows naturally from (\ref{Eqn_v_1}) that $\mathbf{1}_n^T \mathbf{v}^k = \sum_{i=1}^{n} v_i^k =1$ holds for $k \geq K$, which implies $v_i^k \in [\eta^{n-1}/n, 1]$ for $i=1,\cdots, n$ and $k \geq K$.

Given $\mathbf{\bar{P}}^{k} = (\mathbf{\bar{V}}^{k+1})^{-1} \mathbf{\bar{C}}^{k} \mathbf{\bar{V}}^{k}$, we have $p_{ij}^k=c_{ij}^k{v_j^k}/{v_i^{k+1}}$. Further using $v_i^k \in [\eta^{n-1}/n, 1]$, we have $p_{ij}^k \geq \eta^n/n$ for $j \in \mathcal{N}_I^{out} \cup \{i\}$ and $p_{ij}^k = 0$ for $j \notin \mathcal{N}_I^{out} \cup \{i\}$.

Denote $\mathbf{P}_{t:k}$ as $\mathbf{P}_{t:k}=\mathbf{P}^{t} \mathbf{P}^{t-1} \cdots \mathbf{P}^{k}$, one can obtain
\begin{equation*}
	\begin{aligned}
		& [(\mathbf{I}_n - \mathbf{1}_n (\mathbf{v}^{k+1})^T) \otimes \mathbf{I}_d] \mathbf{P}_{k:k-N_P+1} \\
		& = [(\mathbf{\bar{P}}_{k:k-N_P+1} - \mathbf{1}_n(\mathbf{v}^{k-N_P+1})^T) (\mathbf{I}_n - \mathbf{1}_n (\mathbf{v}^{k-N_P+1})^T) ] \\
		& \quad \otimes \mathbf{I}_d
	\end{aligned}
\end{equation*}
Therefore, following Lemma 4 in \cite{nedic2009distributed}, we have
\begin{equation*}
	\begin{aligned}
		& \|[(\mathbf{I}_n - \mathbf{1}_n (\mathbf{v}^{k+1})^T) \otimes \mathbf{I}_d] \mathbf{a}\| \\
		& \leq \|(\mathbf{\bar{P}}_{k:k-N_P+1} - \mathbf{1}_n(\mathbf{v}^{k-N_P+1})^T) \otimes \mathbf{I}_d \| \\ 
		& \quad \ \cdot \| [(\mathbf{I}_n - \mathbf{1}_n (\mathbf{v}^{k-N_P+1})^T) \otimes \mathbf{I}_d] \mathbf{b} \| \\
		& \leq r_P \|[(\mathbf{I}_n - \mathbf{1}_n (\mathbf{v}^{k-N_P+1})^T) \otimes \mathbf{I}_d ] \mathbf{b}\|
	\end{aligned}
\end{equation*}
for any $k \geq K+N_P-1$ and any vector $\mathbf{b} \in \mathbb{R}^{nd}$. \hfill{$\blacksquare$}

\begin{Lemma}\label{Lemma_R}
	Under Assumptions \ref{Assumption_directed_graph} and \ref{Assumption_objective_functions}, and the parameter design in Table II, the sequence of row-stochastic matrices $\{\mathbf{\bar{R}}^k\}$ for $k \geq K$ is ergodic, and $\{\boldsymbol{\phi}^k\}$ is the unique absolute probability sequence for $\{\mathbf{\bar{R}}^k\}$. Denote $Q_R$ as $Q_R=2n(1+\eta^{-(n-1)})/(1-\eta^{n-1})$ and define a positive integer $N_R$ such that $r_{R}\triangleq Q_R (1-\eta^{n-1})^{\frac{N_R-1}{n-1}} <1$ holds. Then for any $k\geq K+N_R-1$ and any vector $\mathbf{b}\in \mathbb{R}^{nd}$, vector $\mathbf{a} \triangleq \mathbf{R}^{k} \mathbf{R}^{k-1} \cdots \mathbf{R}^{k-N_R+1}\mathbf{b}$ always satisfies $\|[(\mathbf{I}_n - \mathbf{1}_n (\boldsymbol{\phi}^{k+1})^T) \otimes \mathbf{I}_d] \mathbf{a}\| \leq r_R \|[(\mathbf{I}_n - \mathbf{1}_n (\boldsymbol{\phi}^{k-N_R+1})^T) \otimes \mathbf{I}_d ] \mathbf{b}\|$.
\end{Lemma}

{\it Proof}: Following a similar line of reasoning for Lemma \ref{Lemma_P}, Lemma \ref{Lemma_R} can be easily obtained and hence we omit the proof here.\hfill{$\blacksquare$}

\section*{Appendix B: Proof of Lemma \ref{Lemma_x}}

{\it Proof}: From (\ref{system-matrix}), we have $\mathbf{x}^{k+1} = \mathbf{R}_{k:k-\bar{N}+1} \mathbf{x}^{k-\bar{N}+1} - \lambda (\mathbf{A}^k\mathbf{y}^k + \sum_{l=1}^{\bar{N}-1}\mathbf{R}_{k:k-l+1}\mathbf{A}^{k-l}\mathbf{y}^{k-l} )$ for $k\geq K+\bar{N}-1$ where $\mathbf{R}_{k:k-\bar{N}+1} = \mathbf{R}^{k} \mathbf{R}^{k-1} \cdots \mathbf{R}^{k-\bar{N}+1}$. Then, it follows
\begin{equation}\label{Eqn_x_1}
	\begin{aligned}
		& \|\mathbf{\tilde{x}}_{\text{w}}^{k+1}\| \leq \big\| \big[ \big(\mathbf{I}_n - \mathbf{1}_n (\boldsymbol{\phi}^{k+1})^T\big) \otimes \mathbf{I}_d \big] \mathbf{R}_{k:k-\bar{N}+1} \mathbf{x}^{k-\bar{N}+1} \big\| \\
		& \ + \lambda \big\| \big[ \big(\mathbf{I}_n - \mathbf{1}_n (\boldsymbol{\phi}^{k+1})^T\big) \otimes \mathbf{I}_d \big] \mathbf{A}^k\mathbf{y}^k \big\| \\
		& \ + \lambda \sum_{l=1}^{\bar{N}-1}\big\| \big[ \big(\mathbf{I}_n - \mathbf{1}_n (\boldsymbol{\phi}^{k+1})^T\big) \otimes \mathbf{I}_d \big] \mathbf{R}_{k:k-l+1} \mathbf{A}^{k-l} \mathbf{y}^{k-l} \big\| \\
		& \leq r_R \| \mathbf{\tilde{x}}_{\text{w}}^{k-\bar{N}+1} \| + \lambda Q_R \sqrt{n} {\textstyle\sum_{l=0}^{\bar{N}-1}} \|\mathbf{y}^{k-l}\|
	\end{aligned}
\end{equation}
where in the derivation we used Lemma \ref{Lemma_R} and the facts $\| \mathbf{I}_n - \mathbf{1}_n (\boldsymbol{\phi}^{k+1})^T\| \leq 2\sqrt{n}$, $\|\mathbf{A}^k\| \leq \sqrt{n}$, and $2\sqrt{n}\leq Q_R$. Following Lemma 1 in \cite{saadatniaki2020decentralized}, we can obtain
\begin{equation}\label{Eqn_y}
	\begin{aligned}
		\| \mathbf{y}^{k}\| \leq n \bar{\beta} \|\mathbf{\tilde{x}}_{\text{w}}^k\| + n \bar{\beta} \|\mathbf{r}^k\| + \|\mathbf{\tilde{s}}_{\text{w}}^k \| 
	\end{aligned}
\end{equation}
for $k\geq K$. Further combing (\ref{Eqn_x_1}) and (\ref{Eqn_y}) leads to the result in (\ref{Eqn_lemma_x}) for $k\geq K+\bar{N}-1$. \hfill{$\blacksquare$}

\section*{Appendix C: Proof of Lemma \ref{Lemma_r}}

{\it Proof}: Given $\mathbf{r}^k = \mathbf{1}_n \otimes \mathbf{\tilde{x}}_{\text{w}}^{k} - \mathbf{1}_n \otimes \mathbf{x}^* = \big[\mathbf{1}_n (\boldsymbol{\phi}^k)^T \otimes \mathbf{I}_d \big] \mathbf{x}^{k} - \mathbf{1}_n \otimes \mathbf{x}^*$, for $k\geq K$ we have
\begin{equation}\label{Eqn_r_1}
	\begin{aligned}
		& \|\mathbf{r}^{k+1}\| \leq \big\|\big[\mathbf{1}_n (\boldsymbol{\phi}^{k+1})^T \otimes \mathbf{I}_d \big] \mathbf{R}^k \mathbf{x}^k - \mathbf{1}_n \otimes \mathbf{x}^* \\
		& \qquad - \lambda \big[(\mathbf{1}_n (\boldsymbol{\phi}^{k+1})^T\mathbf{\bar{A}}^k \mathbf{v}^k \mathbf{1}_n^T) \otimes \mathbf{I}_d \big] \mathbf{y}^k \big\| \\
		& \qquad + \lambda \big\|\mathbf{1}_n (\boldsymbol{\phi}^k)^T \otimes \mathbf{I}_d\big\| \cdot \|\mathbf{A}^k\| \cdot \big\|\mathbf{y}^k -(\mathbf{v}^k \mathbf{1}_n^T \otimes \mathbf{I}_d) \mathbf{y}^k\big\| \\
		& \leq \big\| \big[\mathbf{1}_n (\boldsymbol{\phi}^{k})^T \otimes \mathbf{I}_d \big] \mathbf{x}^k - \mathbf{1}_n \otimes \mathbf{x}^* - \lambda \delta^k (\mathbf{1}_{n \times n} \otimes \mathbf{I}_d) \mathbf{y}^k \big\| \\
		& \qquad + \lambda n \|\mathbf{\tilde{s}}_{\text{w}}^{k} \|\\
	\end{aligned}
\end{equation}
where $\delta^k=(\boldsymbol{\phi}^{k+1})^T\mathbf{\bar{A}}^k \mathbf{v}^k$, and in the derivation we used the facts $\|\mathbf{1}_n (\boldsymbol{\phi}^k)^T \otimes \mathbf{I}_d\| \leq \sqrt{n}$, $\|\mathbf{A}^k\| \leq \sqrt{n}$, $\mathbf{y}^k -(\mathbf{v}^k \mathbf{1}_n^T \otimes \mathbf{I}_d) \mathbf{y}^k = (\mathbf{\bar{V}}^k \otimes \mathbf{I}_d) \mathbf{\tilde{s}}_{\text{w}}^{k}$, and $\|\mathbf{\bar{V}}^k \otimes \mathbf{I}_d \| \leq 1$.

Now we show that $\delta^k \in [\eta^{n-1}/n, \, 1]$. From the proof of Lemma \ref{Lemma_P} in Appendix A, we have $v_i^k \in [\eta^{n-1}/n, \, 1]$ for $k\geq K$ where $v_i^k$ denotes the $i$-th entry of vector $\mathbf{v}^k$. Given $\sum_{j=1}^{n}a_{ij}^k=1$ for $k \geq K$, one can obtain ${\eta^{n-1}}/{n} \leq \sum_{j=1}^{n}a_{ij}^k v_j^k \leq 1$ for each $i=1,\cdots, n$. Therefore, we have $\delta^k=\sum_{i=1}^{n} \phi_i^k \sum_{j=1}^{n}a_{ij}^k v_j^k \in [{\eta^{n-1}}/{n}, \, 1]$ based on the fact $\sum_{i=1}^{n} \phi_i^k=1$.

Next we focus on the second last term on the right hand side of (\ref{Eqn_r_1}), which can be expressed as follows
\begin{equation}\label{Eqn_r_2}
	\begin{aligned}
		& \big\| \big[\mathbf{1}_n (\boldsymbol{\phi}^{k})^T \otimes \mathbf{I}_d \big] \mathbf{x}^k - \mathbf{1}_n \otimes \mathbf{x}^* - \lambda \delta^k (\mathbf{1}_{n \times n} \otimes \mathbf{I}_d) \mathbf{y}^k \big\| \\
		& \leq \big\| \mathbf{1}_n \otimes \mathbf{\bar{x}}_{\text{w}}^{k} - \mathbf{1}_n \otimes \mathbf{x}^* - \lambda \delta^k (\mathbf{1}_{n \times n} \otimes \mathbf{I}_d) \nabla f(\mathbf{1}_n \otimes \mathbf{\bar{x}}_{\text{w}}^{k}) \big\| \\
		& \quad + \lambda \delta^k \big\| (\mathbf{1}_{n \times n} \otimes \mathbf{I}_d) (\nabla f^k - \nabla f(\mathbf{1}_n \otimes \mathbf{\bar{x}}_{\text{w}}^{k})) \big\| \\
		& \leq \sqrt{n} \big\|\mathbf{\bar{x}}_{\text{w}}^{k} - \mathbf{x}^* - \lambda \delta^k \nabla F(\mathbf{\bar{x}}_{\text{w}}^{k}) \big\| + \lambda n \bar{\beta} \| \mathbf{\tilde{x}}_{\text{w}}^{k}\|
	\end{aligned}
\end{equation}
As the global objective function $F$ is $\alpha_F$-strongly convex and $\beta_F$-smooth with $\alpha_F\leq \beta_F$, from Lemma 10 in \cite{qu2017harnessing}, we obtain
\begin{equation}\label{Eqn_r_3}
	\begin{aligned}
		\big\|\mathbf{\bar{x}}_{\text{w}}^{k} - \mathbf{x}^* - \lambda \delta^k \nabla F(\mathbf{\bar{x}}_{\text{w}}^{k}) \big\| \leq \frac{1}{\sqrt{n}}\big(1- \lambda n^{-1} \eta^{n-1}\alpha_F\big) \|\mathbf{r}^k\|
	\end{aligned}
\end{equation}
if $\lambda \leq 1/ \beta_F$ holds. Therefore, combining (\ref{Eqn_r_1}), (\ref{Eqn_r_2}), and (\ref{Eqn_r_3}), we have that the inequality in (\ref{Eqn_lemma_r}) holds for $k\geq K$ if $\lambda$ satisfies $\lambda \leq {1}/\beta_F$. \hfill{$\blacksquare$}

\section*{Appendix D: Proof of Lemma \ref{Lemma_s}}

{\it Proof}: Denoting $\mathbf{z}^k$ as $\mathbf{z}^k = \nabla f^{k+1}-\nabla f^{k}$, from (\ref{Eqn_system_matrix_form_3}) we have $\mathbf{s}^{k+1} = \mathbf{P}_{k:k-\bar{N}+1} \mathbf{s}^{k-\bar{N}+1} + ((\mathbf{\bar{V}}^k)^{-1}\otimes \mathbf{I}_d) \mathbf{B}^k \mathbf{z}^k + \sum_{l=1}^{\bar{N}-1} \mathbf{P}_{k:k-l+1} ((\mathbf{\bar{V}}^{k-l+1})^{-1}\otimes \mathbf{I}_d) \mathbf{B}^{k-l} \mathbf{z}^{k-l}$ for $k\geq K+ \bar{N}-1$. Given $\mathbf{\tilde{s}}_{\text{w}}^k = \mathbf{s}^k - (\mathbf{1}_n(\mathbf{v}^k)^T \otimes \mathbf{I}_d) \mathbf{s}^k$, one can obtain the following inequality for $k\geq K+ \bar{N}-1$
\begin{equation}\label{Eqn_s_1}
	\begin{aligned}
		& \|\mathbf{\tilde{s}}_{\text{w}}^{k+1}\| \leq r_P \big\| \big[ \big(\mathbf{I}_n - \mathbf{1}_n (\mathbf{v}^{k-\bar{N}+1})^T\big) \otimes\mathbf{I}_d \big] \mathbf{s}^{k-\bar{N}+1} \big\| \\
		& \qquad + \big\| \mathbf{I}_n - \mathbf{1}_n (\mathbf{v}^{k+1})^T\big\| \cdot \big\|( \mathbf{\bar{V}}^k)^{-1} \big\| \cdot \big\| \mathbf{B}^k \big\| \cdot \big\| \mathbf{z}^k \big\| \\
		& \qquad + {\textstyle \sum_{l=1}^{\bar{N}-1}} Q_P \big\| ( \mathbf{\bar{V}}^{k-l+1})^{-1} \big\| \cdot \big\| \mathbf{B}^{k-l} \big\| \cdot \big\| \mathbf{z}^{k-l} \big\| \\
		& \leq r_P \|\mathbf{\tilde{s}}_{\text{w}}^{k-\bar{N}+1}\| + (n\sqrt{n}Q_P/\eta^{n-1}) {\textstyle \sum_{l=0}^{\bar{N}-1}} \| \mathbf{z}^{k-l}\|
	\end{aligned}
\end{equation}
where we used the facts $\big\| \mathbf{I}_n - \mathbf{1}_n (\mathbf{v}^{k+1})^T\big\| \leq 2\sqrt{n} \leq Q_P$, $\big\|( \mathbf{\bar{V}}^k)^{-1} \big\| \leq {n}/{\eta^{n-1}}$, $\big\| \mathbf{B}^k \big\| \leq \sqrt{n}$, and Lemma \ref{Lemma_P} in Appendix A.

Using the relationship $\mathbf{z}^k = \nabla f^{k+1}-\nabla f^{k}$, we can obtain the following relationship for $k \geq K$
\begin{equation}\label{Eqn_s_2}
	\begin{aligned}
		& \|\mathbf{z}^k\| \leq \big( {\textstyle\sum_{i=1}^{n}} \beta_i^2 \big\|\mathbf{x}_i^{k+1} - \mathbf{x}_i^{k} \big\|^2 \big)^{1/2} \ \leq \ \bar{\beta} \big\|\mathbf{x}^{k+1} - \mathbf{x}^{k} \big\|\\
		& \leq \bar{\beta} \big\|(\mathbf{R}^{k}-\mathbf{I}) \mathbf{x}^{k}\big\|+ \lambda \bar{\beta} \big\|\mathbf{A}^k \mathbf{y}^k \big\|\\
		& \leq \bar{\beta} \big\|(\mathbf{R}^{k}-\mathbf{I}) \big[(\mathbf{I}_n - \mathbf{1}_n (\boldsymbol{\phi}^k)^T ) \otimes \mathbf{I}_d \big] \mathbf{x}^{k}\big\|+ \lambda \bar{\beta} \big\|\mathbf{A}^k \mathbf{y}^k \big\|\\
		& \leq (2\sqrt{n}\bar{\beta} + \lambda n\sqrt{n}\bar{\beta}^2) \|\mathbf{\tilde{x}}_{\text{w}}^k \| + \lambda n \sqrt{n}\bar{\beta}^2 \|\mathbf{r}^k\| + \lambda \sqrt{n}\bar{\beta} \| \mathbf{\tilde{s}}_{\text{w}}^k \|
	\end{aligned}
\end{equation}
where we used the facts $\big\|\mathbf{R}^k -\mathbf{I} \big\| \leq 2\sqrt{n}$, $\big((\mathbf{I}_n - \mathbf{1}_n (\boldsymbol{\phi}^k)^T ) \otimes \mathbf{I}_d \big) \mathbf{x}^{k}= \mathbf{\tilde{x}}_{\text{w}}^k$, and (\ref{Eqn_y}) in the derivation.

Therefore, from (\ref{Eqn_s_1}) and (\ref{Eqn_s_2}), we can obtain the inequality in (\ref{Eqn_lemma_s}) for $k\geq K+ \bar{N}-1$. \hfill{$\blacksquare$}

\bibliographystyle{plain}
\bibliography{abbr_bibli}

\end{document}